\documentclass[preprint,12pt]{elsarticle}
\usepackage{graphicx}
\usepackage{subfigure}
\usepackage{epstopdf}
\usepackage{rotating}
\usepackage{framed,multirow}
\usepackage{amsmath}
\usepackage{amssymb}
\usepackage{color}
\usepackage{algorithm}
\usepackage{algorithmic}
\usepackage{algorithm}
\usepackage{algorithmic}

\numberwithin{equation}{section}

\graphicspath{{figure/}}
\begin{document}
\title{Enforcing exact boundary and initial conditions in the deep mixed residual method\footnote{The first two authors contributed equally to the work.}}	

	\author[label1,label2]{Liyao Lyu}
	\ead{lyuliyao@msu.edu}
	\author[label1]{Keke Wu}
	\ead{wukekever@gmail.com}
	\author[label1,label4]{Rui Du\corref{cor1}}
	\ead{durui@suda.edu.cn}
	\author[label1,label4]{Jingrun Chen\corref{cor1}}
	\ead{jingrunchen@suda.edu.cn}	
	\address[label1]{School of Mathematical Sciences, Soochow University, Suzhou, 215006, China}
	\address[label2]{CW Chu College, Soochow University, Suzhou, 215006, China}
	\address[label4]{Mathematical Center for Interdisciplinary Research, Soochow University, Suzhou, 215006, China}	
	\cortext[cor1]{Corresponding authors.}

%\subjclass[2000]{Primary 54C40, 14E20; Secondary 46E25, 20C20}

%\date{January 1, 2001 and, in revised form, June 22, 2001.}


%\keywords{Differential geometry, algebraic geometry}

\begin{abstract}
	In theory, boundary and initial conditions are important for the wellposedness of partial differential equations (PDEs). Numerically, these conditions can be enforced exactly in classical numerical methods, such as finite difference method and finite element method. Recent years have witnessed growing interests in solving PDEs by deep neural networks (DNNs), especially in the high-dimensional case. However, in the generic situation, a careful literature review shows that boundary conditions cannot be enforced exactly for DNNs, which inevitably leads to a modeling error. In this work, based on the recently developed deep mixed residual method (MIM), we demonstrate how to make DNNs satisfy boundary and initial conditions automatically in a systematic manner. As a consequence, the loss function in MIM is free of the penalty term and does not have any modeling error. Using numerous examples, including Dirichlet, Neumann, mixed, Robin, and periodic boundary conditions for elliptic equations, and initial conditions for parabolic and hyperbolic equations, we show that enforcing exact boundary and initial conditions not only provides a better approximate solution but also facilitates the training process.
\end{abstract}
    \maketitle

%%====================== Section 1 Introduction ========================
\section{Introduction}
	Partial differential equation (PDE) is one of the most important tools to model various phenomena in science, engineering, and finance. It has been a long history of developing reliable and efficient numerical methods for PDEs. Notable examples include finite difference method~\cite{leveque_finite_2007}, finite element method~\cite{zienkiewicz1977finite}, and discontinuous Galerkin method~\cite{cockburn2012discontinuous}. For low-dimensional PDEs, these methods are proved to be accurate and demonstrated to be efficient. However, they run into the curse of dimensionality for high-dimensional PDEs, such as Schr\"odinger equation in the quantum many-body problem \cite{dirac1981principles}, Hamilton-Jacobi-Bellman equation in the stochastic optimal control \cite{bardi2008optimal}, and nonlinear Black-Scholes equation for pricing financial derivatives \cite{hull2009options}. 
	
	In the last decade, significant advancements in deep learning have driven the development of solving PDEs in the framework of deep learning, especially in the high-dimensional case where deep neural networks overcome the curse of dimensionality by construction; see~\cite{E2017Dec,E2018Mar,han2018solving,sirignano2018dgm:,raissi2019physics-informed,chen2019quasi-monte,Beck2019Aug,Hutzenthaler2019Mar,zang2020weak,Becker2020May, Lyu2020Jun} for examples and references therein. Among these, deep Ritz method uses the variational form (if exists) of the corresponding PDE as the loss function~\cite{E2018Mar} and deep Galerkin method (DGM) uses the PDE residual in the least-squares senses as the loss function~\cite{sirignano2018dgm:}. It is worth mentioning that DGM has no connection with Galerkin from the perspective of numerical PDEs although it is named after Galerkin. In~\cite{raissi2019physics-informed}, physics-informed neural networks is proposed to combine observed data with PDE models. The mixed residual method (MIM) first rewrites a PDE into a first-order system and then uses the system residual in the least-squares sense as the loss function~\cite{Lyu2020Jun}. 
	These progresses demonstrate the strong representability of deep neural networks (DNNs) for solving PDEs.
	
	In classical numerical methods, basis functions or discretization stencils have compact supports or sparse structures. Machine-learning methods, instead, employ DNNs as trial functions, which are globally defined. This stark difference makes DNNs overcome the curse of dimensionality while classical numerical methods cannot. However, there are still unclear issues for DNNs, such as the dependence of approximation accuracy on the solution regularity and the enforcement of exact boundary conditions. It is straightforward to enforce exact boundary conditions in classical numerical methods while it is highly nontrivial for DNNs due to their global structures. A general strategy is to add a penalty term in the loss function which penalizes the discrepancy between a DNN evaluated on the boundary and the exact boundary condition. Such a strategy inevitably introduces a modeling error which pollutes the approximation accuracy and typically has a negative impact on the training process \cite{Chen2020May}. Therefore, it is always desirable to construct DNNs which automatically satisfy boundary conditions and there are several efforts towards this objective \cite{5061501,Berg2018Nov,Sheng2020Apr}. It is shown that Dirichlet boundary condition can be enforced exactly over a complex domain in~\cite{Berg2018Nov}. This idea cannot be applied for Neumann boundary condition since the solution value on the boundary is not available. This issue is solved by constructing the trail DNN in a different way~\cite{5061501}. However, for mixed boundary condition, this construction has a serious issue at the intersection of Dirichlet and Neumann boundary conditions and an approximation has to be applied. Therefore, it is so far that an exact enforcement of mixed boundary condition for DNNs has still been lacking.
		
	In this work, in the framework of MIM, we demonstrate how to make DNNs satisfy boundary and initial conditions automatically in a systematic manner. As a consequence, the loss function in MIM is free of penalty term and does not have any modeling error. The success relies heavily on the unique feature of MIM. In MIM, both the PDE solution and its derivatives are treated as independent variables, very much like the discontinuous Galerkin method \cite{cockburn2012discontinuous} and least-squares finite element method \cite{Bochev2015}, while other deep-learning methods only have PDE solutions as unknown variables. Therefore, it is straightforward to enforce Dirichlet boundary condition for all deep-learning methods and Neumann boundary condition requires a bit more efforts. But for mixed boundary condition, it is only possible in MIM to enforce the exact condition since a direct access to both the solution and its first-order derivatives is only available simultaneously in MIM. For completeness, we also study Robin and periodic boundary conditions for elliptic equations. For parabolic equations, the enforcement of exact initial conditions only requires the PDE solution and thus can be done for all machine-learning methods. For wave equations, the enforcement of exact initial conditions requires both the PDE solution and its first-order derivatives with respect to time and thus only can be done in MIM. Using numerous examples, including Dirichlet, Neumann, mixed, Robin, and periodic boundary conditions for elliptic equations, and initial conditions for parabolic and hyperbolic equations, we show that enforcing exact boundary and initial conditions not only provides a better approximate solution but also facilitates the training process. 	
	
	The paper is organized as follows. In Section \ref{sec: MIM}, we will give a complete description of MIM and also DGM for the comparison purpose since both methods work for general types of PDEs. Constructions of DNNs for different boundary conditions are derived in Section \ref{sec:bc} with numerical demonstrations in Section \ref{sec:numerical result}. Constructions of DNNs for initial conditions are derived with numerical demonstrations in Section \ref{sec:numerical result initial}. Conclusions are drawn in Section \ref{sec:conclusion}.

%%====================== Section 2 MIM ======================== 

\section{Deep mixed residual method}
\label{sec: MIM}

	For completeness, we will introduce MIM in this section. We will also introduce DGM for the comparison purpose. In both methods, there are three main components. Interested readers may refer to \cite{Lyu2020Jun} and \cite{sirignano2018dgm:} for details.
	\begin{itemize}
		\item Modeling: Rewrite the original problem into an optimization problem by defining a loss function;
		\item Architecture: Build the trail function space using DNNs;
		\item Optimization: Search for the optimal set of parameters in the DNN which minimizes the loss function.
	\end{itemize}  

	\subsection{Modeling}
	\label{subsection:Modeling}
	
	Consider an elliptic equation as an example
	\begin{equation*}
	\left\{
	\begin{aligned}
	&-\nabla \cdot (a(x)\nabla  u) = f(x) & x\in \Omega\subset\mathbb{R}^d\\
	& \Gamma u(x)= g(x) & x\in \partial \Omega
	\end{aligned}
	\right.
	\end{equation*}
	with different types of boundary conditions
	\begin{equation*}
		\begin{aligned}
		&\text{Dirichlet}& & \Gamma u = u,\\
		&\text{Neumann}& &\Gamma u = a(x)\nabla u \cdot \nu,\\
		&\text{Robin}& & \Gamma u = a(x)\nabla u \cdot \nu + u.
		\end{aligned}
	\end{equation*}
	Here $\nu$ represents the outward unit normal vector on $\partial\Omega$ and $\nu_k$ is its $k$-th component.
	For completeness, we also consider the mixed boundary condition and periodic boundary condition.
	
	In general, DGM uses the following loss function
	\begin{equation}\label{equ:loss function}
	\begin{aligned}
	L(u) = \|\nabla \cdot (a(x)\nabla u ) - f\|^2_{2,\Omega} + \lambda \|\Gamma u-g\|^2_{2,\partial \Omega},
	\end{aligned}
	\end{equation}
	where $\lambda$ is the penalty parameter that balances the two terms in \eqref{equ:loss function}.
	
	In MIM, the PDE is first rewritten into a first-order system of equations, very much like the discontinuous Galerkin method \cite{cockburn2012discontinuous} and least-squares finite element method \cite{Bochev2015}. The loss function is then defined as the residual of the first-order system in the least-squares sense. Details of MIM for different types of high-order PDEs can be found in \cite{Lyu2020Jun}. The loss function for the elliptic equation is defined as 
	\begin{equation*}
		L(u,p) = \|a(x)\nabla u - p\|_{\Omega,2}^2 + \|\nabla \cdot p + f\|_{\Omega,2}^2 + \lambda \|\Gamma u-g\|^2_{2,\partial \Omega}. 
	\end{equation*}
	At the moment, we still use the penalty term to enforce the boundary condition for brevity and will show that MIM can be free of penalty in all cases in Section \ref{sec:bc} and Section \ref{sec:numerical result initial}. Afterwards, multiple DNNs are used to approximate the PDE solution and its derivatives: one DNN $\hat{u}_\theta$ is used to approximate the PDE solution $u$ and the other DNN $\hat{p}_\theta$ is used to approximate its derivatives $\nabla u$. Be aware that $\hat{u}_\theta$ and $\hat{p}_\theta$ are treated as independent variables, which is crucial for the success of enforcing exact boundary and initial conditions in all cases.
	
	\subsection{Architecture}

\begin{figure}[ht]
		\centering
		\includegraphics[width=0.3\textwidth]{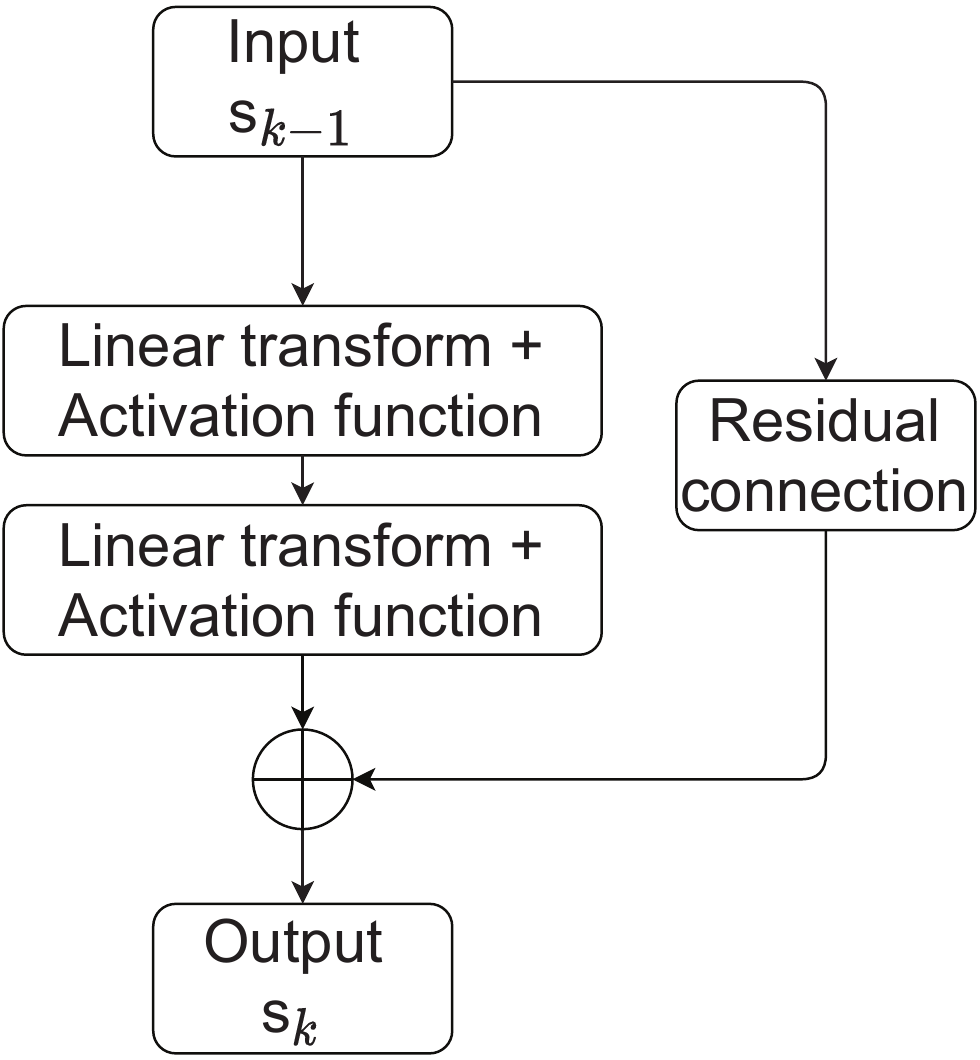}
		\caption{One block of ResNet. A deep neural network contains a sequence of blocks, each of which consists of two fully-connected layers and one shortcut connection.}
		\label{fig:resnet}
	\end{figure}

	ResNet \cite{He2015} is used to approximate the PDE solution and its high-order derivatives. A ResNet consists of $m$ blocks in the following form
	\begin{equation}\label{equ:resnet}
	s_k =\sigma(W_{2,k}\sigma(W_{1,k}s_{k-1}+b_{1,k})+b_{2,k}) +s_{k-1}, \quad k=1,2,\cdots,m.
	\end{equation}
	Here $s_k, b_{1,k}, b_{2,k} \in \mathbb{R}^n$, $W_{1,k}, W_{2,k}\in \mathbb{R}^{n\times n}$. $m$ is the depth of network, $n$ is the width of network, and $\sigma$ is the (scalar) activation function. In numerical tests, we use ReQu $(\max\{x,0\})^2$, ReCu $(\max\{x,0\})^3$ or $swish(x) = x / (1 + \exp(-x))$ as the activation function. The last term on the right-hand side of \eqref{equ:resnet} is called the shortcut connection or residual connection. Each block has two linear transforms, two activation functions, and one shortcut; see Figure \ref{fig:resnet} for demonstration. Such a structure can automatically solve the notorious problem of vanishing/exploding gradient \cite{He2015}.

	Note that the network input is in $\mathbb{R}^d$, $s_i$ is in $\mathbb{R}^n$ and $n$ is often different from $d$, while the output network is in $\mathbb{R}$ for the solution $u$ or in $\mathbb{R}^d$ for $\nabla u$. Therefore, we can pad $x$ by a zero vector or apply a linear transform to get the network input $s_0$ and apply a linear transform from $\mathbb{R}^m$ to $\mathbb{R}$ (or $\mathbb{R}^d$) to make the network output be in the same dimension as the target function $u$ (or $p$). Such a ResNet is denoted as $N_\theta(x)$. All the parameters, including parameters in \eqref{equ:resnet}, in the first linear transform (if exists) and the last linear transform, will be optimized. The total number of parameters in the set is $(2m-1)n^2 + (2m+d+1)n +1$ in DGM and $(4m-2)n^2 + (4m+3d+1)n +d + 1$ in MIM. Note that $d$ shall be redefined as the sum of the spatial dimension and the temporal dimension $1$ for time-dependent problems.

	\subsection{Optimization}
	Take the loss function in MIM of elliptic equation for example
	\begin{equation*}
		\min_{\theta} L(\hat{u}_{\theta},\hat{p}_{\theta}).
	\end{equation*}
	To find the optimal set of parameters in the above problem, we adopt ADAM \cite{kingma2015adam}, a gradient-based method which is widely used for performing machine-learning tasks. Algorithm~\ref{alg:Adam} is provided here for completeness. Note that $g_k^2$ in Algorithm \ref{alg:Adam} denotes the elementwise square of $g_k$, i.e., $g_k\odot g_k$. We will also use $\oplus$ and $\ominus$ to denote the elementwise addition and subtraction in what follows. 
	\begin{algorithm}
	\caption{ADAM}
	\label{alg:Adam}
	\begin{algorithmic}
		\REQUIRE $\alpha$: Learning rate is set to be $0.001$ in this paper.
		\REQUIRE $\beta_1,\beta_2 \in [0,1)$: Exponential decay rates for moment estimates, which are set to be $0.9$ and $0.999$ here, respectively.
		\REQUIRE $\varepsilon \leftarrow 1e-8$
		\REQUIRE $\theta_0$: Initial parameter vector,\\
			\quad$m_0\leftarrow0$ (Initialize the 1-st moment vector),\\
			\quad$v_0\leftarrow0$ (Initialize the 2-nd moment vector),\\
			\quad$k\leftarrow0$ (Initialize the counter)
			
			\WHILE{$\theta_k$ not converged} 
			\STATE $k \leftarrow k+1$ 
			\STATE $x \leftarrow$ random sampling points in $\Omega$ 
			(Update sampling points)
			\STATE $L_k(\theta_{k-1}) \leftarrow L(\hat{u}_{\theta_{k-1}}(x),\hat{p}_{\theta_{k-1}}(x))$ 
			(Evaluate the loss function)
			\STATE $g_k \leftarrow \nabla_\theta L_k(\theta_{k-1})$  
			(Calculate the gradient w.r.t. $\theta$ at step k)
			
			\STATE $m_k \leftarrow \beta_1 \cdot m_{k-1} +(1-\beta_1)\cdot g_k$
			(Update the biased first moment)
			\STATE $v_k \leftarrow \beta_2\cdot v_{k-1} + (1-\beta_2)\cdot g_k^2$
			(Update the biased second raw moment)
			\STATE $\widehat{m}_k \leftarrow m_k/(1-\beta_1^k)$
			(Compute the bias-corrected first moment)
			\STATE $\widehat{v}_k \leftarrow v_k/(1-\beta_2^k)$ 
			(Compute the bias-corrected second raw moment)
			\STATE $\theta_k \leftarrow \theta_{k-1} - \alpha \cdot \widehat{m}_k/(\sqrt{\widehat{v}_k} +\varepsilon)$
			(Update parameters)
			\ENDWHILE
		\RETURN $\theta_k$ (Return parameters)
	\end{algorithmic}
	\end{algorithm}

	\section{Enforcement of exact boundary conditions}\label{sec:bc}
	In this section, we will present technical details of MIM on how to enforce exact boundary conditions as well as available strategies in the literature for comparison.
	
	\subsection{Dirichlet boundary conditon}
	For Dirichlet boundary condition, we have the direct access to the solution value on the boundary. Therefore, it is straightforward to design a DNN that satisfies the exact boundary condition; see~\cite{Berg2018Nov}. Given $u(x) = g(x)$ on $\partial \Omega$, the trial DNN is constructed as
	\begin{equation}
	\label{equ:BCconstruction}			
	\hat{u}_{\theta} = L_D(x)N_{\theta}(x) + G_D(x).
	\end{equation} 
	Here $L_D(x)$ is the distance function to the boundary and thus $L_D(x) = 0, \forall x\in\partial\Omega$, $G_D(x)$ is a (smooth) extension of $g(x)$ in $\Omega$ and $G_D(x)=g(x),\forall x\in\partial\Omega$. It is easy to check that \eqref{equ:BCconstruction} satisfies the boundary condition automatically. 
	
	The choice of $L_D(x)$ is not unique, but there are two necessary conditions that $L_D(x)$ has to satisfy
	\begin{itemize}
		\item  $L_D(x) = 0, \forall x\in\partial\Omega$;
		\item $L_D(x) \neq 0, \forall x\in\Omega$.
	\end{itemize}
    The first condition guarantees the enforcement of Dirichlet boundary condition, and the second condition is important for the approximation accuracy. If there exists $x_0\in\Omega$, such that $L_D(x_0)=0$, then $\hat{u}_{\theta}(x_0) = G_D(x_0)$, which is not necessary to be close to $u(x_0)$.
	
	In DGM, \eqref{equ:BCconstruction} is directly used as the trial function. In MIM, we use \eqref{equ:BCconstruction} for $\hat{u}_{\theta}$ but keep $\hat{p}_{\theta}$ free since $\hat{u}_{\theta}$ and $\hat{p}_{\theta}$ are treated independently.
	
	\subsection{Neumann boundary condition}
	Unfortunately, the construction \eqref{equ:BCconstruction} fails to satisfy other boundary conditions since the direct access to the solution value on the boundary is not available. To illustrate this, we consider Neumann boundary condition of the form $a\nabla u \cdot \nu = g(x) $. According to \eqref{equ:BCconstruction}, the trial DNN solution is constructed as
	\begin{equation*}			
	\hat{u}_{\theta} = L_N(x)N_{\theta}(x) + G_N(x).
	\end{equation*}
	 Similarly, we ask  $a\nabla (L_N(x)N_\theta(x))\cdot \nu = 0$  and $a\nabla G_N(x) \cdot \nu= g(x)$ on $\partial \Omega$. For a general DNN $N_\theta(x)$, the former constraint requires that $L_N(x) =0 $ and $\nabla L_N(x) = 0 $ on $\partial \Omega$. This immediately leads to $\hat{u}_\theta(x) = G_N(x)$ on the boundary. Meanwhile, we can construct $G_N(x)$ such that $a\nabla G_N(x) \cdot \nu= g(x)$ on $\partial \Omega$. However, we do not have the direct access to the solution value on the boundary for Neumann boundary condition. If by choice $G_N(x)\neq u(x), \forall x\in\partial\Omega$, then the discrepancy will be out of control on the boundary and it is impossible to get a good approximation over $\Omega$. 
	 
	 In~\cite{5061501}, a different form of the trail DNN is proposed which satisfies Neumann boundary condition automatically
	 \begin{equation}\label{eqn:nbc}			
	 \hat{u}_{\theta} = L_N(x) F_N(x,N_{\theta}(x)) + N_{\theta}(x).
	 \end{equation}
	The choice of $L_N(x)$ is not unique, but there are two necessary conditions that $L_N(x)$ has to satisfy
	\begin{itemize}
		\item  $L_N(x) = 0, \forall x\in\partial\Omega$;
		\item $a\nabla L_N(x)\cdot\nu \neq 0, \forall x\in\bar{\Omega}=\Omega\cup\partial\Omega$.
	\end{itemize}
	Taking the gradient of \eqref{eqn:nbc} with respect to $x$ yields
	 \begin{equation*}
	 	\nabla \hat{u}_{\theta} =  \nabla L_N(x) F_N(x,N_\theta(x)) + L_N(x)\nabla F_N(x,N_\theta(x)) + \nabla N_\theta(x),
	 \end{equation*}
	where the second term vanishes on $\partial \Omega$. Solving for $F_N(x,N_\theta(x))$ in the above equation on $\partial\Omega$ produces
	 \begin{equation*}
	 	F_N(x,N_\theta(x)) = \frac{g(x) - a\nabla N_\theta(x)\cdot\nu }{a\nabla L_N(x)\cdot\nu }, \quad \forall x\in \partial\Omega,
	 \end{equation*}
	which is extended to the whole domain $\Omega$ as
	 \begin{equation}\label{eqn:nbcG}
		F_N(x,N_\theta(x)) = \frac{G_N(x) - a\nabla N_\theta(x)\cdot\nu }{a\nabla L_N(x)\cdot\nu }, \quad \forall x\in \Omega.
	\end{equation}	 
	 The choice of $L_N(x)$ is important to ensure that the denominator in \eqref{eqn:nbcG} is away from $0$ for any $x\in\bar{\Omega}$. Remember that $G_N(x)$ is a (smooth) extension of $g(x)$ from $\partial\Omega$ to $\Omega$.
	 
	 \eqref{eqn:nbc} and \eqref{eqn:nbcG} will be used in DGM for Neumann boundary condition in Section \ref{sec:numerical result}. In MIM, since we have the direct access to both the solution and its derivative, we can simplify the construction as
     \begin{align}
	         &\hat{p}_\theta = F_N(x,N^*_\theta(x)) \nabla L_N(x)  +  N^*_\theta(x),\label{eqn:nbcmim1}\\
	         & F_N(x,N^*_\theta(x)) = \frac{G_N(x) - a N^*_\theta(x)\cdot\nu}{a \nabla L_N(x)\cdot\nu },\label{eqn:nbcmim2}
	 \end{align}
	 where $N_\theta(x)$ is a one dimensional DNN and $N^*_\theta(x)$ is a $d$ dimensional DNN. Apparently, the construction \eqref{eqn:nbcmim1}-\eqref{eqn:nbcmim2} has more degrees of freedom than \eqref{eqn:nbc}-\eqref{eqn:nbcG}. Thus we expect MIM can provide better approximations than DGM, as will be shown in Section \ref{sec:numerical result}.

	 \subsection{Mixed boundary condition}
	 
	Enforcement of exact mixed boundary condition was also considered in~\cite{5061501}. For mixed boundary condition problem,
		\begin{equation}
		\begin{aligned}
		&u(x) = g_D(x) & x \in \Gamma_D,\\
		&\nabla u(x) \cdot \nu = g_N(x) & x \in \Gamma_N,\\
		\end{aligned}
		\end{equation}
	 the construction starts with the following trail form
	\begin{equation}\label{equ:former construction}			
	\hat{u}_{\theta} = L_D(x)N_\theta(x) + L_D(x)L_N(x) F_N(x,N_{\theta}(x)) + G_D(x),
	\end{equation}
	where $G_D$ is a (smooth) extension of $g_D$. Here $\Gamma_D$ is the portion of $\partial\Omega$ with Dirichlet condition and $\Gamma_N$ is the portion of $\partial\Omega$ with Neumann condition.
	Conditions on $L_N(x)$ and $L_D(x)$ are
	\begin{itemize}
		\item  $L_D(x) =0, \forall x\in\Gamma_D$;
		\item $L_D(x) \neq 0, \forall x\in \bar{\Omega}/\Gamma_D$;
		\item $L_N(x) =0, \forall x\in\Gamma_N$;
		\item $L_N(x) \neq 0, \forall x\in \bar{\Omega}/\Gamma_N$.
	\end{itemize}
	Using these conditions, after some algebraic calculations, we arrive at
	\begin{equation}\label{eqn:mbcG}
	F_N(x,N_\theta(x)) = \frac{ G_N(x) -  L_D(x)\left( a \nabla N_\theta(x)\cdot\nu\right) - N_{\theta}(x)\left(a\nabla L_D(x)\cdot\nu\right)}{L_D(x)\left(a \nabla L_N(x)\cdot\nu \right)},
	\end{equation}
	where $G_N$ is a (smooth) extension of $g_N$.
	However, a serious issue is encountered applying \eqref{equ:former construction}-\eqref{eqn:mbcG} since the denominator in \eqref{eqn:mbcG} converges to $0$ when $x$ gets close to the intersection between boundaries of Dirichlet condition and Neumann condition. It was argued in \cite{5061501} that one could add an additional term in the denominator to avoid the numerical issue. However, in the presence of the additional term, the trial DNN will never satisfy the exact mixed boundary condition. This problem remains open so far.
	
	In the framework of MIM, we provide a construction that satisfies the mixed boundary condition automatically.
	\begin{equation}\label{eqn:mbcmim1}
		\begin{aligned}
		&\hat{u}_{\theta} = L_D(x) N_{\theta}(x) +G_D(x),\\
		&\hat{p}_{\theta} = G(x,N^*_{\theta}(x))\nabla L_N(x) + N^*_{\theta}(x),
		\end{aligned}
	\end{equation}
	where 
	\begin{equation}\label{eqn:mbcmin2}
	F_N(x,N^*_\theta(x)) = \frac{G_N(x) - a N^*_{\theta}(x)\cdot\nu}{ a\nabla L_N(x)\cdot\nu}
	\end{equation}
	It is easy to check that \eqref{eqn:mbcmim1}-\eqref{eqn:mbcmin2} satisfies the mixed boundary condition by construction and it will be shown that \eqref{eqn:mbcmim1}-\eqref{eqn:mbcmin2} performs well numerically.
    
    \subsection{Robin boundary condition}
	Consider Robin boundary condition
	\begin{equation}\label{eqn:rbc}
	    a\nabla u \cdot \nu + u = g.
	\end{equation}
	The trial DNN is constructed as 
	\begin{equation}\label{eqn:rbc1}
	    \hat{u}_\theta(x) = L_R(x) F_R(x,N_\theta(x)) + N_\theta(x),
	\end{equation}
	where
	\begin{equation}\label{eqn:rbc2}
	     F_R(x,N_\theta(x)) = \frac{ G_R(x) - N_\theta(x) - a\nabla N_\theta(x) \cdot\nu}{a\nabla L_R(x) \cdot\nu}.
	\end{equation}
	Here $G_R(x)$ is a (smooth) extension of $g(x)$ over $\Omega$.
	The choice of $L_R(x)$ is not unique, but there are two necessary conditions that $L_R(x)$ has to satisfy
	\begin{itemize}
		\item  $L_R(x) = 0, \forall x\in\partial\Omega$;
		\item $a\nabla L_R(x)\cdot\nu \neq 0, \forall x\in\bar{\Omega}=\Omega\cup\partial\Omega$.
	\end{itemize}
	It is not difficult to check that \eqref{eqn:rbc1}-\eqref{eqn:rbc2} satisfies \eqref{eqn:rbc} by taking the gradient of $\hat{u}_\theta$ in \eqref{eqn:rbc1} together with \eqref{eqn:rbc2}. The construction \eqref{eqn:rbc1}-\eqref{eqn:rbc2} will be used in DGM. 
	
	In MIM, we can construct the trail solution as
	\begin{equation}\label{eqn:rbcmin1}
		\begin{aligned}
			    \hat{u}_\theta(x) &= N_\theta(x),\\
			    \hat{p}_\theta(x) &= L_R(x) F_R(x,N^*_\theta(x)) + N^*_\theta(x),
		\end{aligned}
	\end{equation} 
	where
	\begin{equation}\label{eqn:rbcmin2}
	F_N(x,N^*_\theta(x)) = \frac{G_R(x) - N_\theta(x) - a N^*_\theta(x)\cdot \nu}{a\nabla L_R(x)\cdot \nu}.
	\end{equation}
	Similarly, we can check that \eqref{eqn:rbcmin1}-\eqref{eqn:rbcmin2} satisfies the boundary condition \eqref{eqn:rbc}.
	
	\subsection{Periodic boundary condition}
	The construction for periodic boundary condition follows mainly on \cite{han2020solving}. The details are as follows.
	Consider the periodic boundary condition of the form
	\begin{equation}\label{eqn:pbc}
	u(x_1,\cdots,x_{i} + I_i,\cdots,x_d) = u(x_1,\cdots,x_{i},\cdots,x_d),\quad i=1,\cdots,d,
	\end{equation}
	where $I_i$ is the period along the $i$-th direction.
	
	To make DNNs satisfy the periodicity automatically, we construct a transform $T: \mathbb{R}^d \to  \mathbb{R}^{2kd}$ for the input $x = (x_1,\cdots,x_{i},\cdots,x_d)$ before applying the first fully connected layer of the neural network. The component $x_i$ in $x$ is transformed as follows
	\begin{equation*}
	x_i  \to \{\sin( 2\pi j\frac {x_i} { I_i}) , \cos(2\pi j\frac {x_i} { I_i})\}_{j = 1}^k
	\end{equation*}
	for $i = 1, \cdots, d$. This treatment is similar to building the Fourier series of a function, which can capture both high-frequency and low-frequency information of the function.
	\begin{figure}[H]
		\centering
		\includegraphics[width = 1.0 \textwidth]{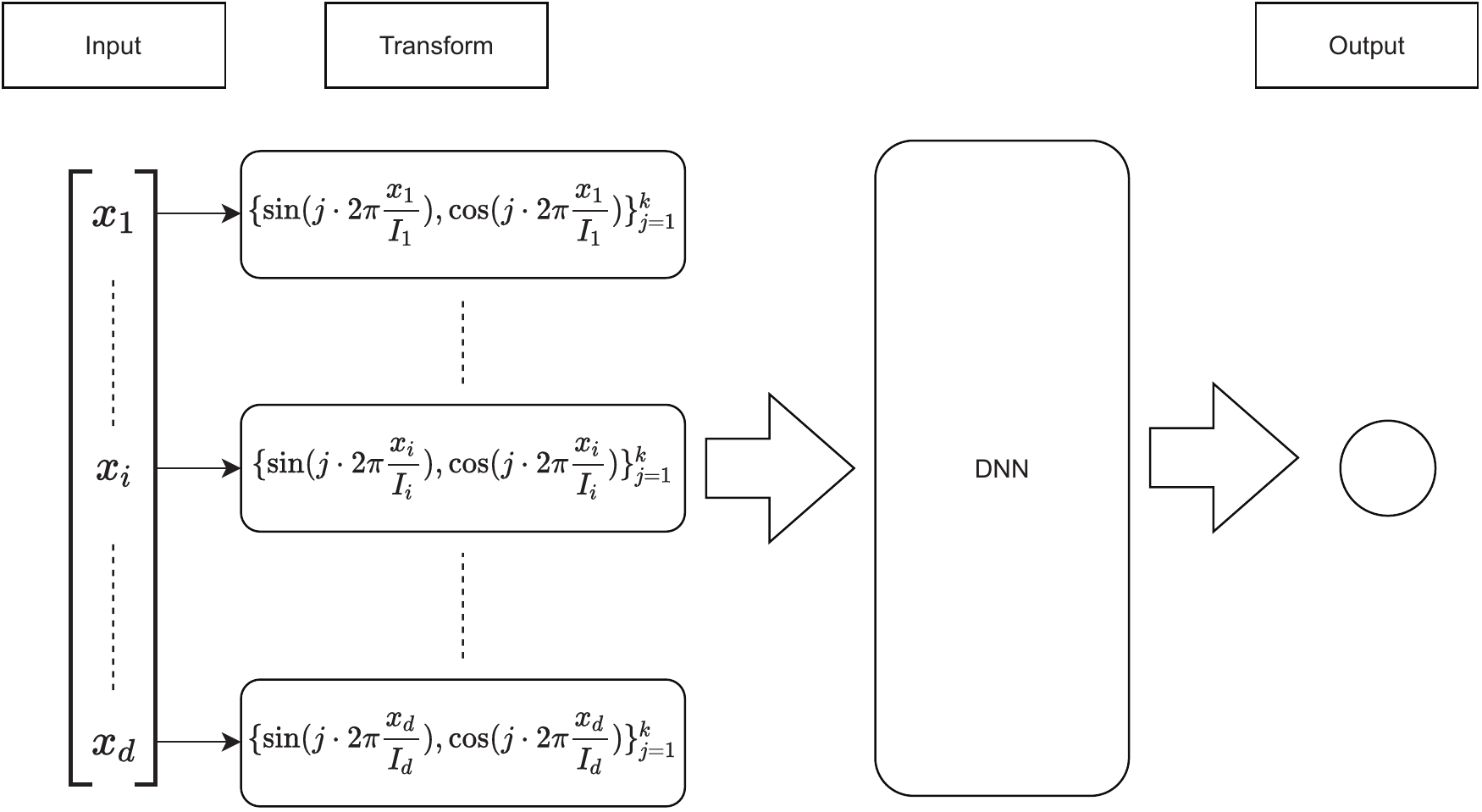}
		\caption{Network structure for periodic boundary condition.}
		\label{fig:periodic}
	\end{figure}
	
%%====================================================================%%
	\section{Numerical results for boundary conditions}
	\label{sec:numerical result}
	
	In this section, we will test the constructions in Section \ref{sec:bc} using a series of examples. For quantitative comparison, we use the relative $L^2$ error defined as 
	\begin{equation}\label{eqn:error}
		\epsilon = \frac{\|\hat{u}_{\theta}-u\|_{2}}{\|u\|_{2}}.
	\end{equation}
	
	\subsection{Dirichlet boundary condition}
	\label{subsec:dirichlet boundary}
	
\begin{figure}
	\centering
	\subfigure{
		\includegraphics[width=0.45\textwidth]{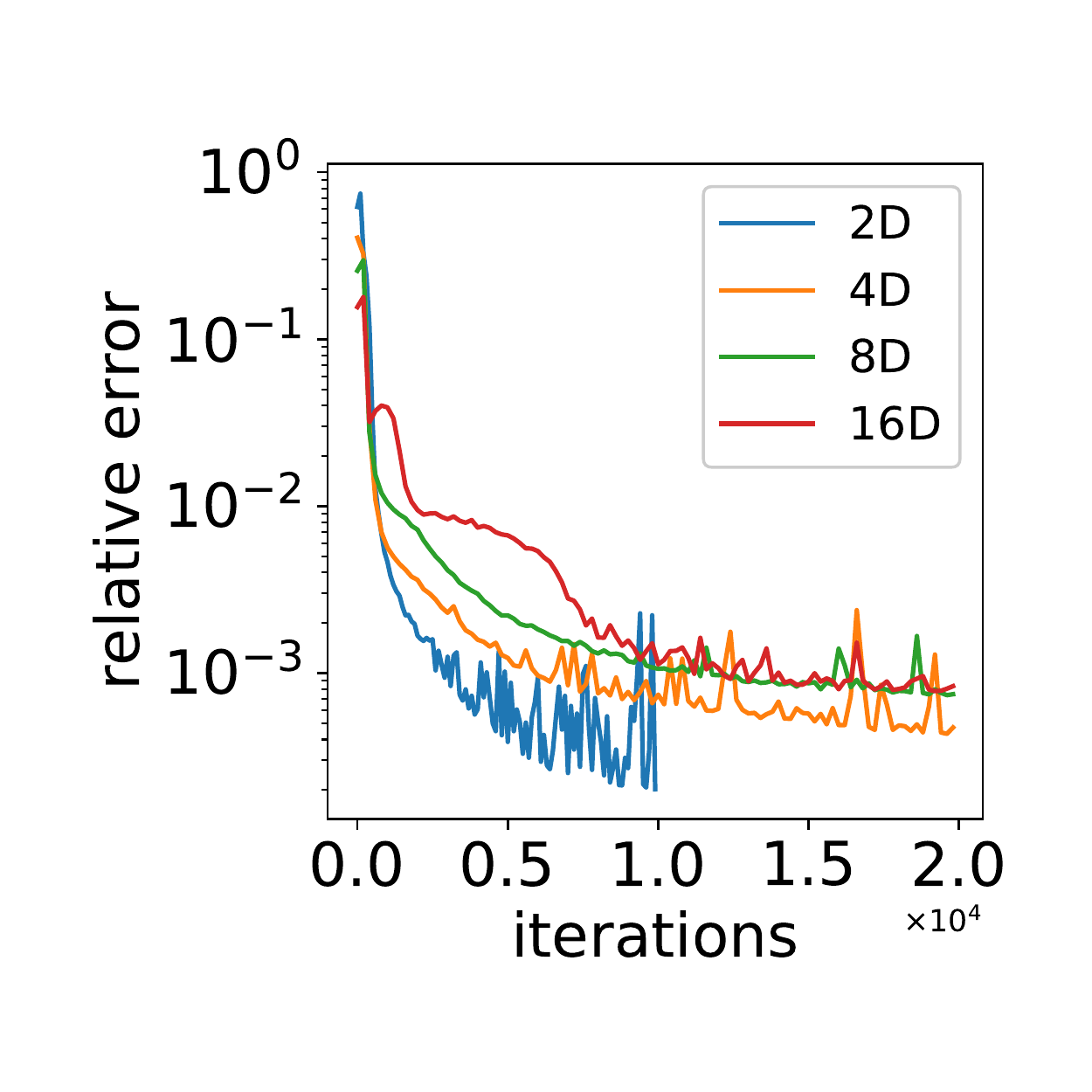}
	}
	\subfigure{
		\includegraphics[width=0.45\textwidth]{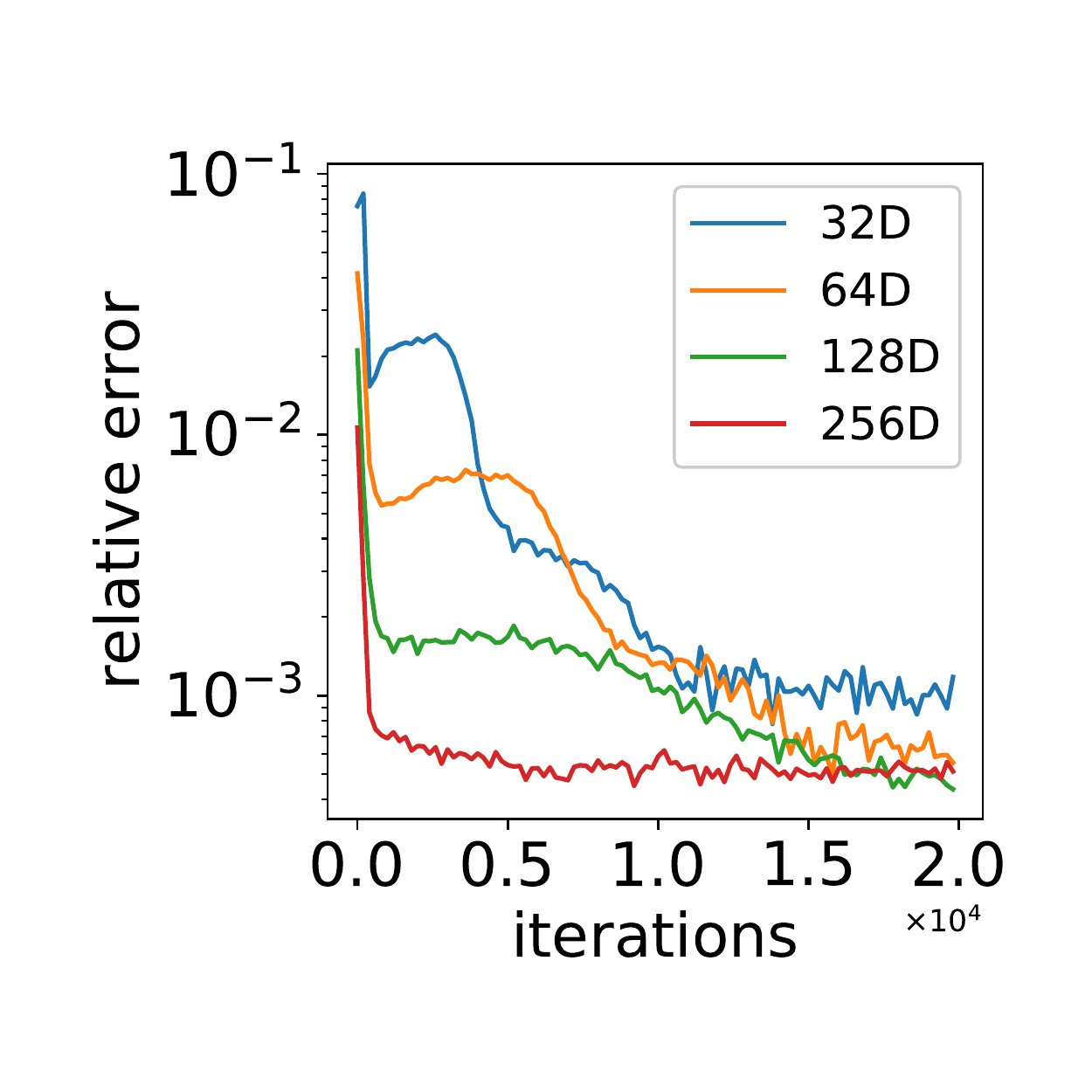}
	}
	\caption{Training processes of MIM for \eqref{eqn:elliptic} with Dirichlet boundary condition in different dimensions. Excellent performance of MIM in high dimensions is observed while the training process of DGM is found to be difficult to converge in this case.}
	\label{fig:poisson Dirichlet}
	\end{figure}

	Consider a nonlinear elliptic equation
	\begin{equation}\label{eqn:elliptic}
	\left\{\begin{aligned}
	&-\Delta u + u^2 = f  & x\in \Omega\\
	&u = e & x\in \partial \Omega
	\end{aligned} \right.
	\end{equation}	
	with exact solution $u(x)= e^{\|x\|^2}$ defined in a sphere $\Omega = \{x:\|x\|<1\}$. For this problem, we set $L_D(x)= \|x\|-1$ and $G(x) = e$ in \eqref{equ:BCconstruction} for the trail DNN. One can verify that it satisfies the boundary condition automatically. Results of DGM and MIN in relative $L^2$ errors are shown in Table \ref{tab:poisson Dirichlet} with training processes in Figure \ref{fig:poisson Dirichlet}.
	\begin{table}[ht]
	\centering
	\caption{Relative $L^2$ errors of DGM and MIM for \eqref{eqn:elliptic} with Dirichlet boundary condition in different dimensions. ReQu is used as the activation function here and $10000$ sampling points are used in domain $\Omega$. The training process ends after $10000-20000$ epochs.}
	\begin{tabular}{|c|c|c|c|c|}
		\hline
		\multirow{2}*{$d$} & \multirow{2}*{$n$} & \multirow{2}*{$m$}  & \multicolumn{2}{c|}{$\epsilon$} \\
		\cline{4-5}
		~& ~ & ~  & MIM & DGM\\
		\hline
		2 & 10 & 2 & 2.37 e-04 & 3.26 e-04\\
		4 & 15 & 2 & 5.85 e-04 & 3.13 e-04\\
		8 & 20 & 2 & 8.10 e-04 & 3.22 e-04\\
		16& 20 & 2 & 8.63 e-04 & 2.31 e-04\\
		32& 35 & 2 & 1.01 e-03 & 1.53 e-04\\
		64& 70 & 2 & 5.85 e-04 & 9.41 e-05\\
		128& 144 & 2 & 4.63 e-04 &  - \\
		256& 280 & 2 & 5.19 e-04 &  - \\				
		\hline
	\end{tabular}
	\label{tab:poisson Dirichlet}
	\end{table}
	Since both MIM and DGM are implemented without the penalty term, the errors are quite small ($<0.1\%$). If a penalty term is used for the boundary condition, then larger errors are observed; see examples in \cite{chen2019quasi-monte} for details. In this example, DGM has slightly better results than MIM overall. However, when the dimension becomes larger, it is found that MIM converges well while DGM fails to converge for a given number of epochs, as shown in the last two lines of Table \ref{tab:poisson Dirichlet}.
	
\begin{figure}
		\centering
		\subfigure[2D]{
			\includegraphics[width=0.45\textwidth]{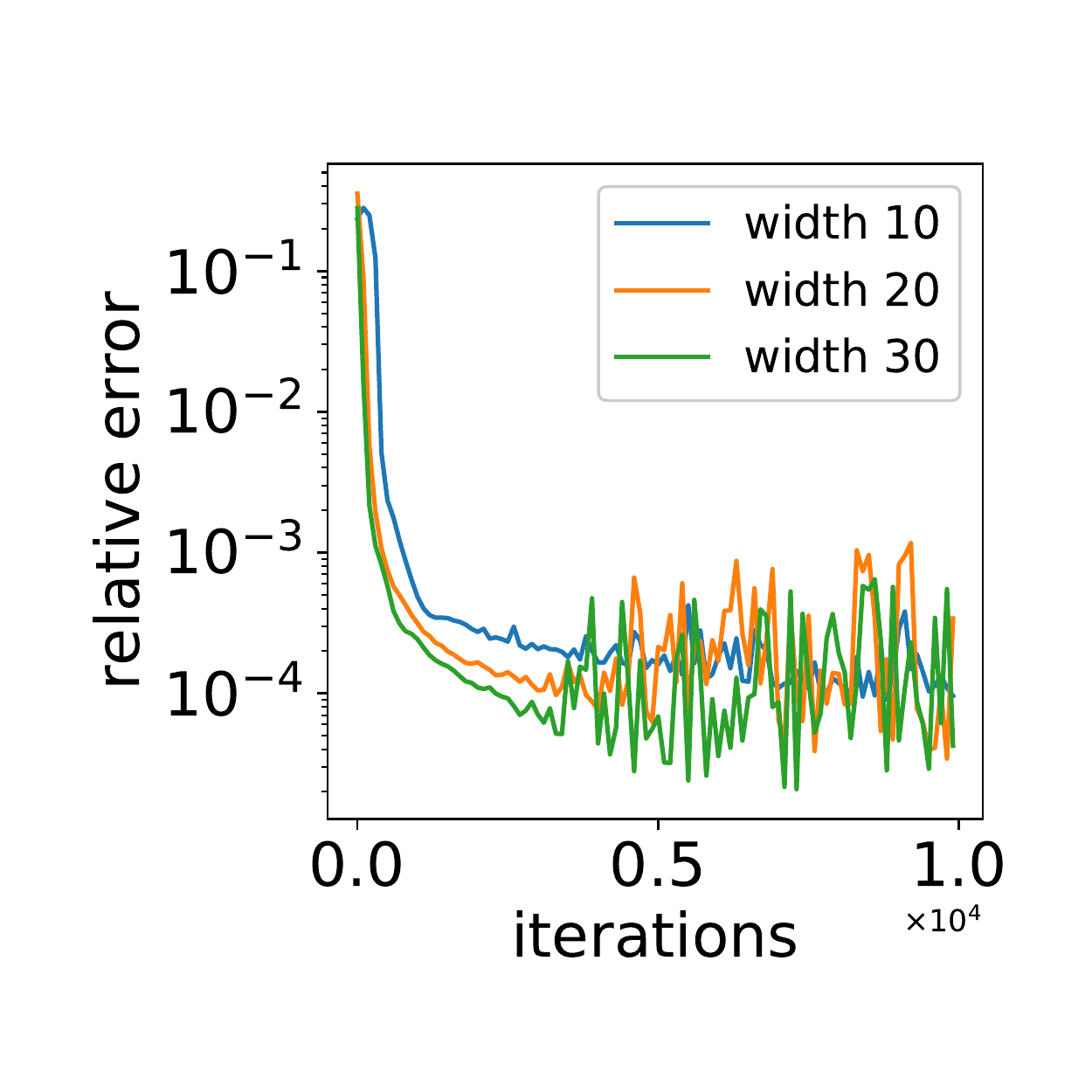}
		}
		\subfigure[4D]{
			\includegraphics[width=0.45\textwidth]{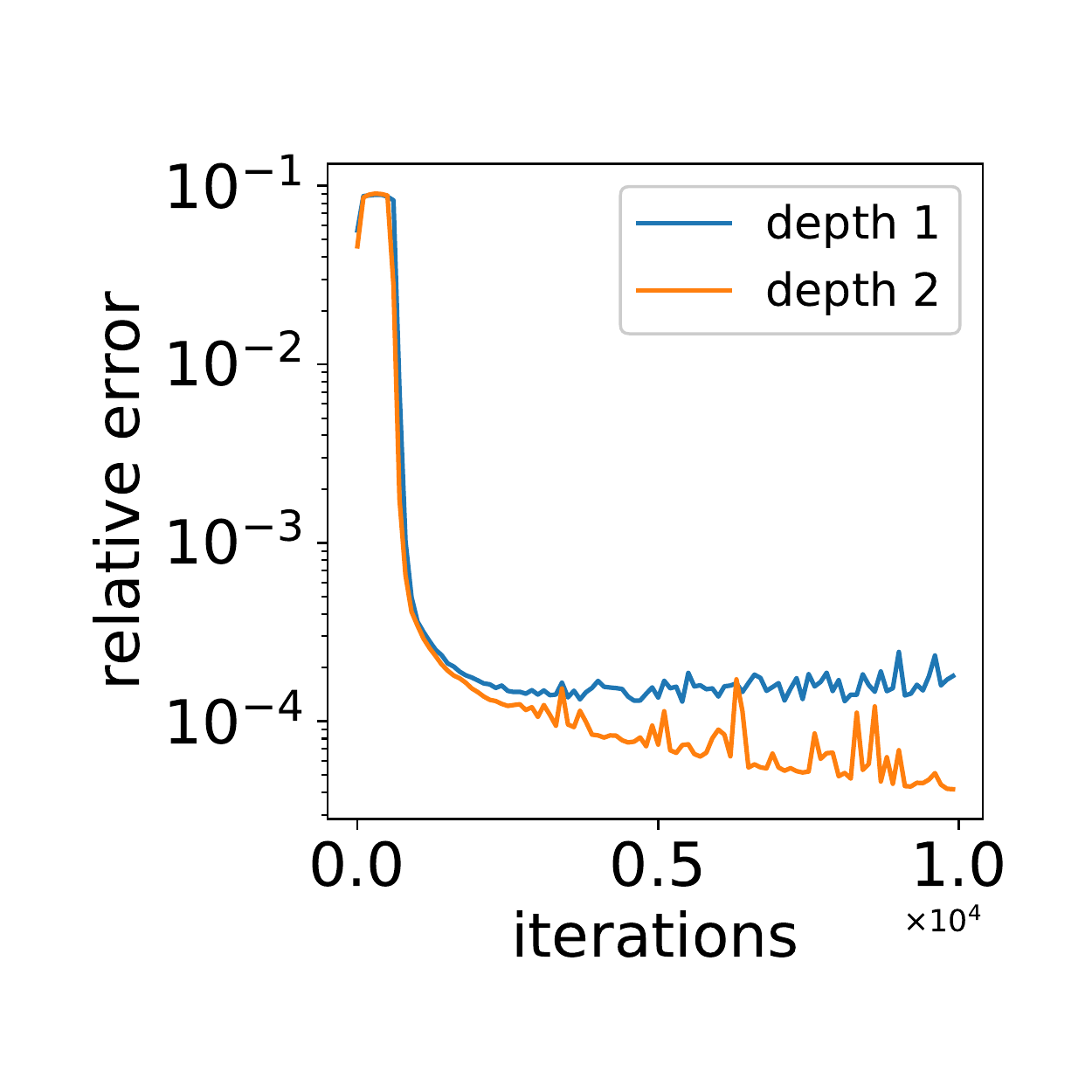}
		}
		\caption{Training processes of MIM for Monge-Amp\'{e}re equation \eqref{equ:MA} with Dirichlet boundary condition in different dimensions.}
		\label{fig:MA Dirichlet}
	\end{figure}

	Next we consider the Monge-Amp\'{e}re equation
	\begin{equation}\label{equ:MA}
	\left\{\begin{aligned}
	&\det(\nabla^2 u) = f(x) & x\in \Omega \\
	&u(x) = e^{\frac{1}{d}} & x\in \partial \Omega
	\end{aligned} \right.
	\end{equation}
	with exact solution $u(x)= e^{\frac{\sum_{i=1}^d x_i^2}{d}}$ over $\Omega = \{x:\|x\|<1\}$. 
	The trial solution in MIM is constructed as 
	\begin{equation*}
	\begin{aligned}
	&\hat{u}_{\theta} = (1-|x|^2)N_\theta (x)+ e^{\frac{1}{d}},\\
	&\hat{p}_{\theta} = N^*_\theta(x).
	\end{aligned}
	\end{equation*}
	Here $\hat{u}_{\theta}$ satisfies the boundary condition automatically and $\hat{p}_{\theta}$ is free. The loss function in MIM is defined as 
	\begin{equation*}
	L(u,p) = \|p - \nabla u \|_{2,\Omega}^2  + \|\det(\nabla p) - f\|_{2,\Omega}^2.
	\end{equation*}
	Results of MIM are shown in Table \ref{tab:MA Dirichlet} with training processes in Figure \ref{fig:MA Dirichlet}. Again, without the penalty term, the errors are quite small ($<0.1\%$). 
	\begin{table}[ht]
		\caption{Relative $L^2$ errors for Monge-Amp\'{e}re equation \eqref{equ:MA} with Dirichlet boundary condition in different dimensions. ReQu is used as the activation function here and $50000$ sampling points are used in domain $\Omega$. The training process ends after $10000$ epochs.}
		\centering
		\begin{tabular}{|c|c|c|c|}
			\hline
			$d$ & $n$ & $m$ & $\epsilon$ \\
			\hline
			2 & 10 & 2 & 1.39 e-04 \\
			2 & 20 & 2 & 2.16 e-04 \\
			2 & 30 & 2 & 1.91 e-04 \\
			4 & 20 & 1 & 1.66 e-04 \\
			4 & 20 & 2 & 6.82 e-05 \\
			\hline
		\end{tabular}
		\label{tab:MA Dirichlet}
	\end{table}
	\subsection{Neumann boundary condition}
	
\begin{figure}
		\centering
		\subfigure[MIM]{
			\includegraphics[width=0.45\textwidth]{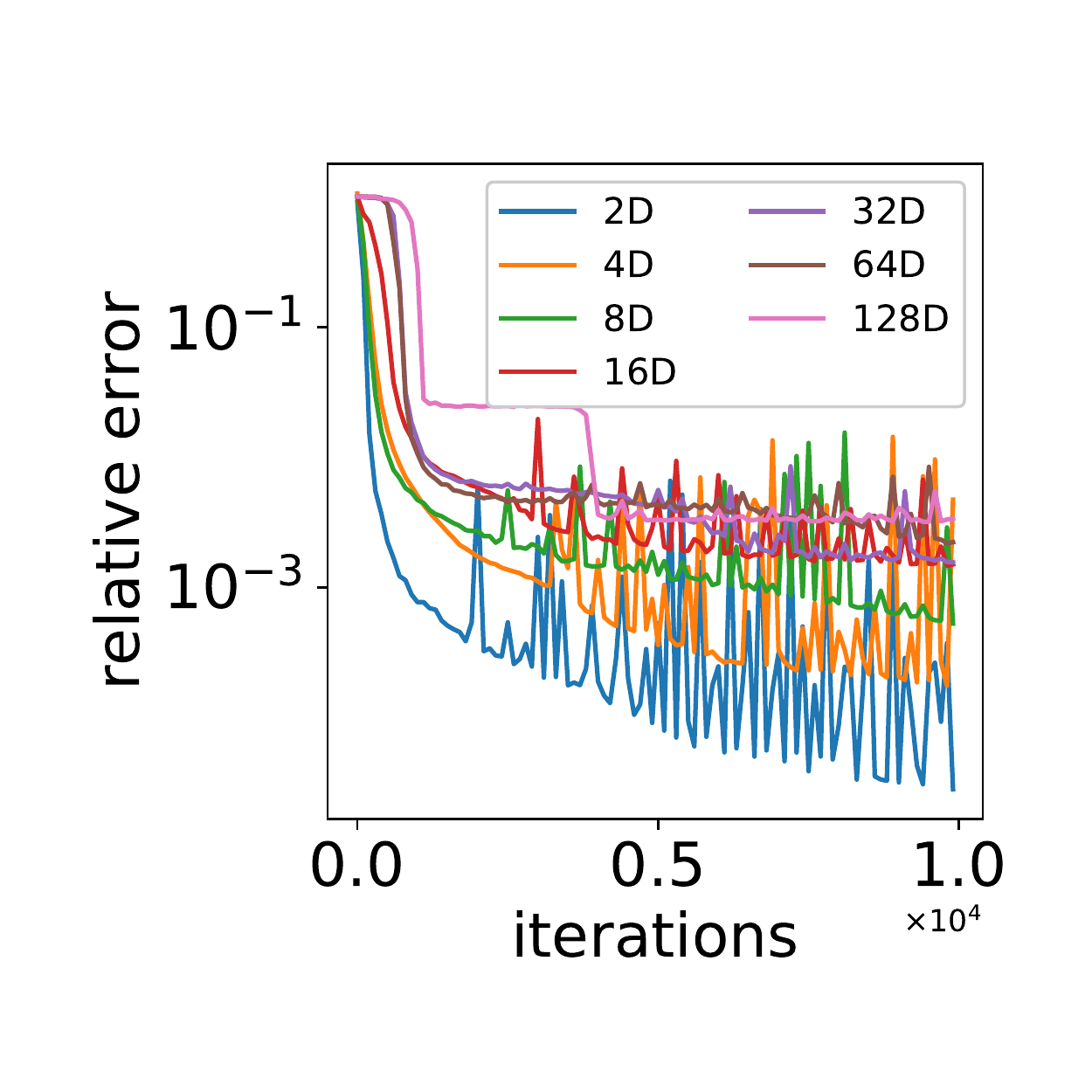}
		}
		\subfigure[DGM]{
			\includegraphics[width=0.45\textwidth]{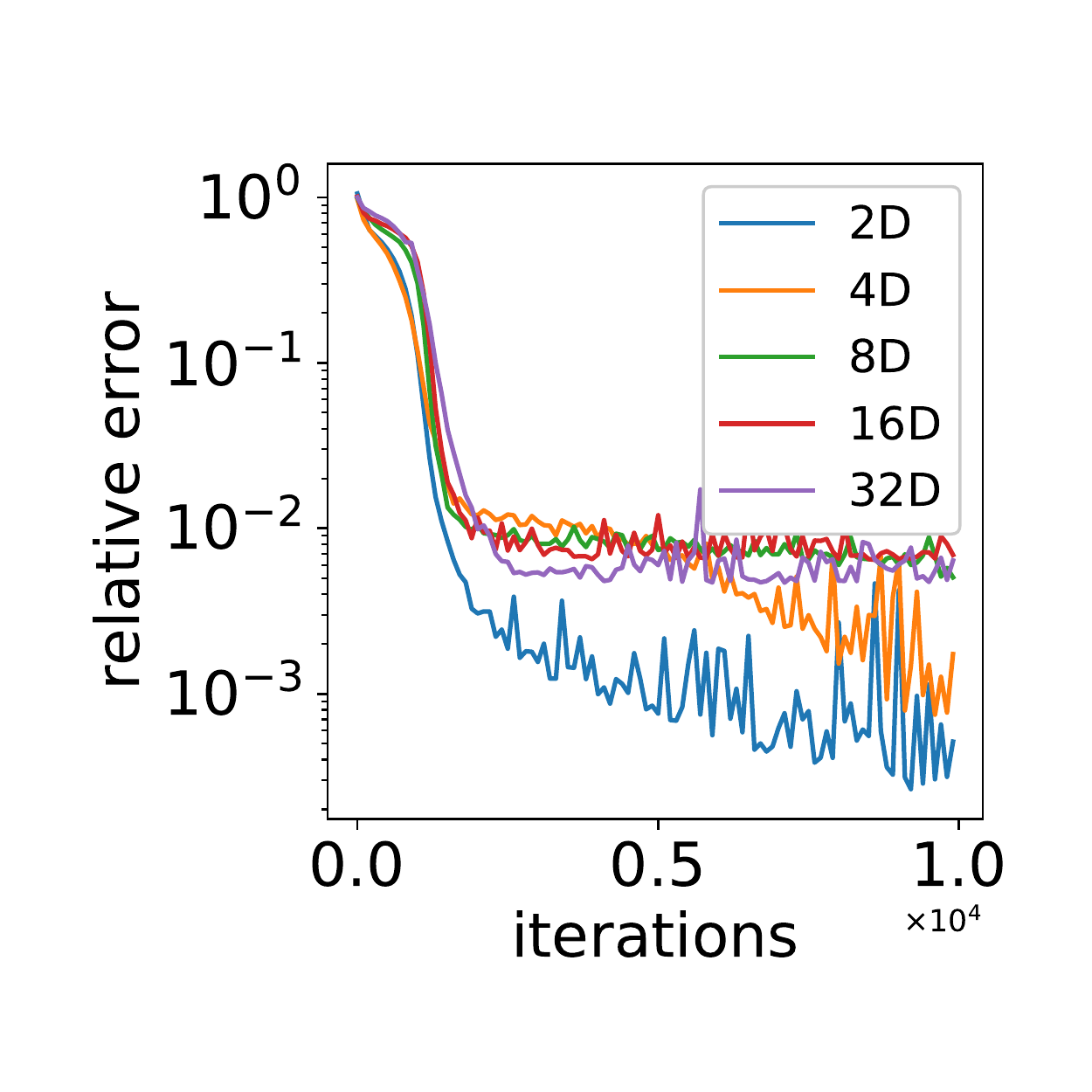}
		}
		\caption{Training processes of MIM and DGM for \eqref{eqn:elliptic nbc} with Neumann boundary condition in different dimensions. DGM fails to converges when $d=64$ and $128$.}
		\label{fig:Neumann boundary condition}
	\end{figure}

	Consider the elliptic problem
	\begin{equation}\label{eqn:elliptic nbc}
	\left\{\begin{aligned}
	& -\Delta u + u = f & x\in \Omega\\
	& \frac{\partial u}{\partial n } = g & x \in \partial \Omega\\ 
	\end{aligned} \right.	
	\end{equation}
	with exact solution $u(x) = \sum_{k=1}^d \exp(x)$ over $\Omega = [0,1]^d$. The trail solution in MIM is constructed as
	\begin{equation*}
	\begin{aligned}
	&\hat{u}_{\theta} = N_\theta (x),\\
	&\hat{p}_{\theta} = (1-|x|^2)N^*_\theta(x) + (e-1)x + 1.
	\end{aligned}
	\end{equation*}
	Here $\hat{p}_{\theta}$ satisfies the boundary condition automatically and $\hat{u}_{\theta}$ is free. The loss function is defined as
	\begin{equation*}
	\begin{aligned}
	L(u,p) = \|\nabla u - p \|_{2,\Omega}^2  + \|\nabla \cdot p - u + f \|_{2,\Omega}^2.
	\end{aligned}
	\end{equation*}
	Relative $L^2$ errors of MIM and DGM are recorded in Table \ref{tbl:Neumann boundary condition} with training processes in Figure \ref{fig:Neumann boundary condition}. Here DGM is implemented with the penalty method and the penalty parameter $\lambda=1$. MIM provides better results than DGM in all dimensions. One may argue that DGM can produce better results by a fine tuning of the penalty parameter. However, this requires a significant work which violates the main purpose of this work. Moreover, the absence of the penalty term actually facilitates the training process, as shown in the last two lines of Table \ref{tbl:Neumann boundary condition}. We find that it is difficult to train the parameters in DGM when $d=64$ and $128$. Similar observations are provided in \cite{Chen2020May}.
	\begin{table}[ht]
		\centering 
		\caption{Relative $L^2$ errors for \eqref{eqn:elliptic nbc} with Neumann boundary condition in different dimensions. ReQu is used as the activation function here. We use $10000-50000$ sampling points in $\Omega$ and $1000$ sampling points along each side of $\partial\Omega$ in DGM. The training process ends after $10000$ epochs.}
		\begin{tabular}{|c|c|c|c|c|}
			\hline
			\multirow{2}*{$d$} & \multirow{2}*{$n$} & \multirow{2}*{$m$}  & \multicolumn{2}{c|}{$\epsilon$} \\
			\cline{4-5}
			~& ~ & ~  & MIM & DGM\\
			\hline
			2  & 10 & 2 & 2.86 e-05 & 3.67 e-04\\
			4  & 15 & 2 & 6.23 e-04 & 1.37 e-03\\
			8  & 20 & 2 & 1.70 e-03 & 6.12 e-03\\
			16 & 25 & 2 & 2.55 e-03 & 7.18 e-03\\
			32 & 35 & 2 & 3.08 e-03 & 6.14 e-03\\
			64 & 70 & 2 & 2.43 e-03 & -\\
			128& 130& 2 & 3.61 e-03 & -\\
			\hline
		\end{tabular}
		\label{tbl:Neumann boundary condition}
	\end{table}
	Consider another Neumann boundary value problem
	\begin{equation}\label{eqn:ellptic nbc2}
	\left\{
	\begin{aligned}
	&-\Delta u - u = f 
	& x\in \Omega \\
	&\frac{\partial u}{\partial n} = 0 
	&x \in \partial \Omega \\
	\end{aligned}\right.
	\end{equation}
	with exact solution $ \cos(\sum_{k=1}^d x_k^2 -1)$ over $\Omega = \{x:\|x\|<1\}$.

\begin{figure}
	\centering
	\subfigure[MIM]{
		\includegraphics[width=0.45\textwidth]{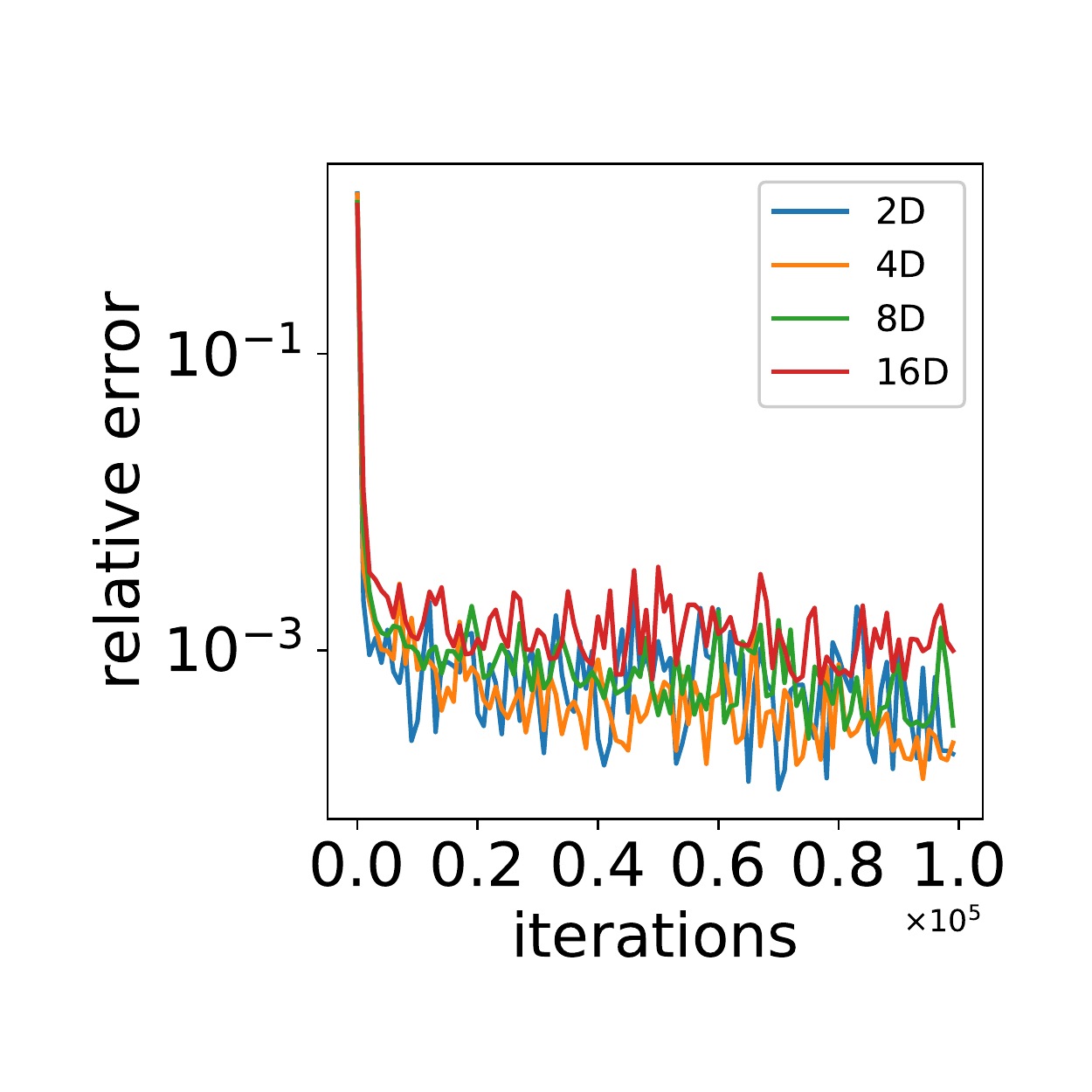}
	}
	\subfigure[DGM]{
		\includegraphics[width=0.45\textwidth]{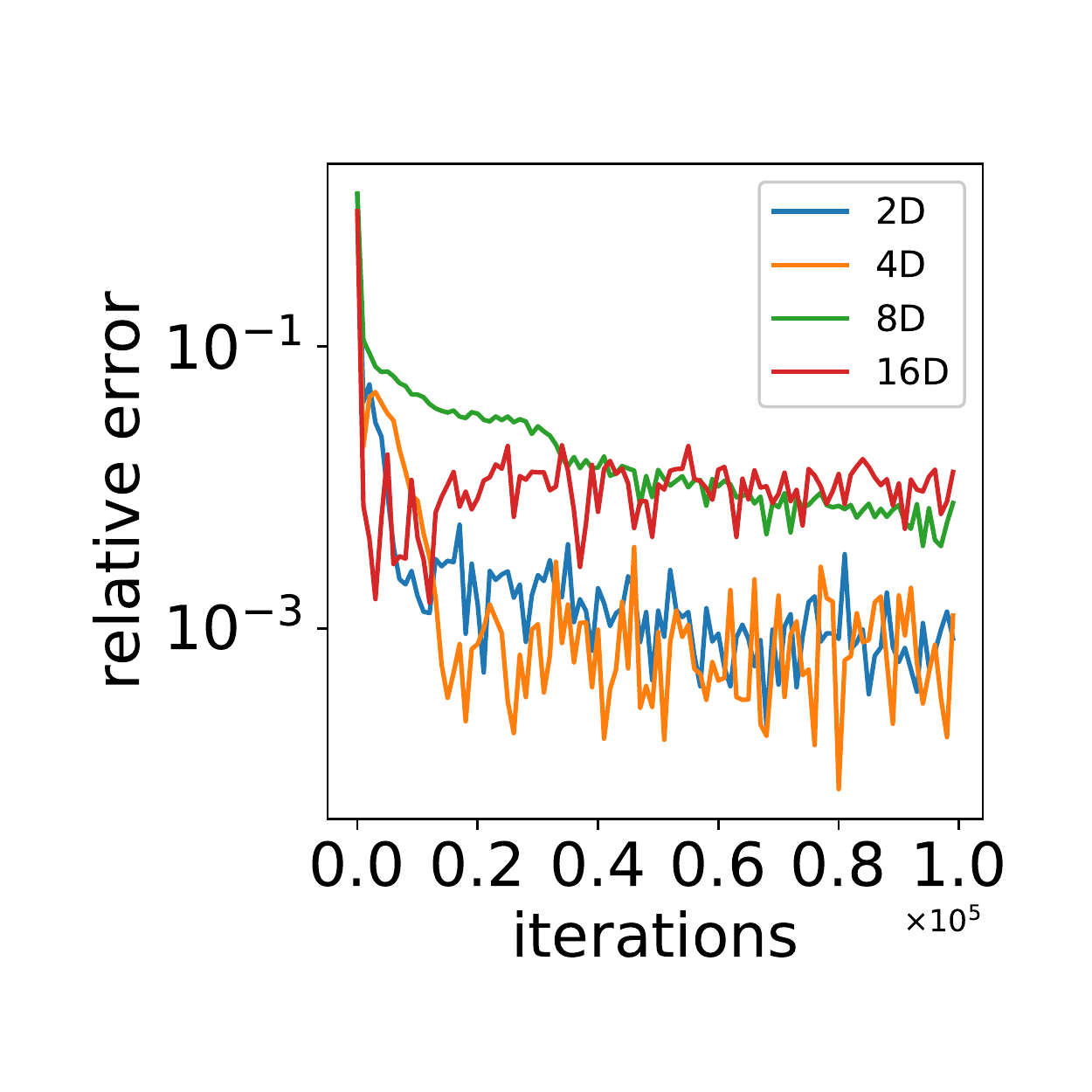}
	}
	\caption{Training processes of MIM and DGM for \eqref{eqn:ellptic nbc2} with Neumann boundary condition in different dimensions.}
	\label{fig:Neumann error general}
\end{figure}

	To remove the penalty term in DGM, following \eqref{equ:former construction}, we construct the trail solution in DGM  as 
	\begin{equation*}
	\hat{u}_{\theta} = \frac{1}{2}\left(\sum_{k=1}^d x_k^2 -1\right)\left(- \sum_{k=1}^d x_k \frac{\partial N_\theta(x)}{\partial x_k}\right) + N_{\theta}(x).
	\end{equation*}
	Consequently, the loss function in DGM is defined as 
	\begin{equation*}
	L(u) = \|\Delta u + u+f \|^2_{2, \Omega},
	\end{equation*}
	which is free of the penalty term. In MIM, $\hat{u}_{\theta} = N_{\theta}(x)$ is free and 
	\begin{equation*}
	\hat{p}_{\theta} = x\left(- x\cdot N^*_{\theta}(x) \right) + N^*_{\theta}(x).
	\end{equation*}
	Results of both methods are shown in Table \ref{tbl:Neumann error general} with training processes in Figure \ref{fig:Neumann error general}. Both methods are now free of the penalty term, but still MIM outperforms DGM in all dimensions. MIM has the direct access to both the solution and its derivatives, and thus provides better approximations.
	\begin{table}
		\centering
		\caption{Relative $L^2$ errors for \eqref{eqn:ellptic nbc2} with Neumann boundary condition in different dimensions. ReQu is used as the activation function here and $10000$ sampling points are used in $\Omega$. The training process ends after $100000$ epochs.}
		\begin{tabular}{|c|c|c|c|c|}
			\hline
			\multirow{2}*{$d$} & \multirow{2}*{$n$} & \multirow{2}*{$m$}  & \multicolumn{2}{c|}{$\epsilon$} \\
			\cline{4-5}
			~& ~ & ~  & MIM & DGM\\
			\hline
			2   & 10 & 3 & 5.19 e-04 & 1.01 e-03\\
			4   & 15 & 3 & 3.60 e-04 & 6.53 e-04\\
			8   & 20 & 3 & 5.84 e-04 & 6.00 e-03\\
			16  & 25 & 3 & 1.14 e-03 & 9.97 e-03\\
			\hline
		\end{tabular}
		\label{tbl:Neumann error general}
	\end{table}
	\subsection{Robin boundary condition}

	Consider the elliptic equation
	\begin{equation}
	\label{equ:Robin1}\left\{
	\begin{aligned}
	&-\Delta u + \pi^2 u = f \quad x\in \Omega\\
	&\frac{\partial u}{\partial \nu} +u = g \quad x\in \partial \Omega
	\end{aligned}\right.
	\end{equation}
	with exact solution $u(x)= \sin(\sum_{k=1}^d x_k) $ over $\Omega=[0,1]^d$.
	
	The constructions \eqref{eqn:rbc}-\eqref{eqn:rbcmin2} provide a systematic way to enforce the exact Robin boundary condition. Since MIM has the direct access to both the solution and its derivatives, we will demonstrate two different ways to construct trail DNNs that satisfy the Robin boundary condition automatically.
	
\begin{figure}
	\centering
	\subfigure[MIM \eqref{eqn:rbc new1}-\eqref{eqn:rbc new3}]{
		\includegraphics[width=0.45\textwidth]{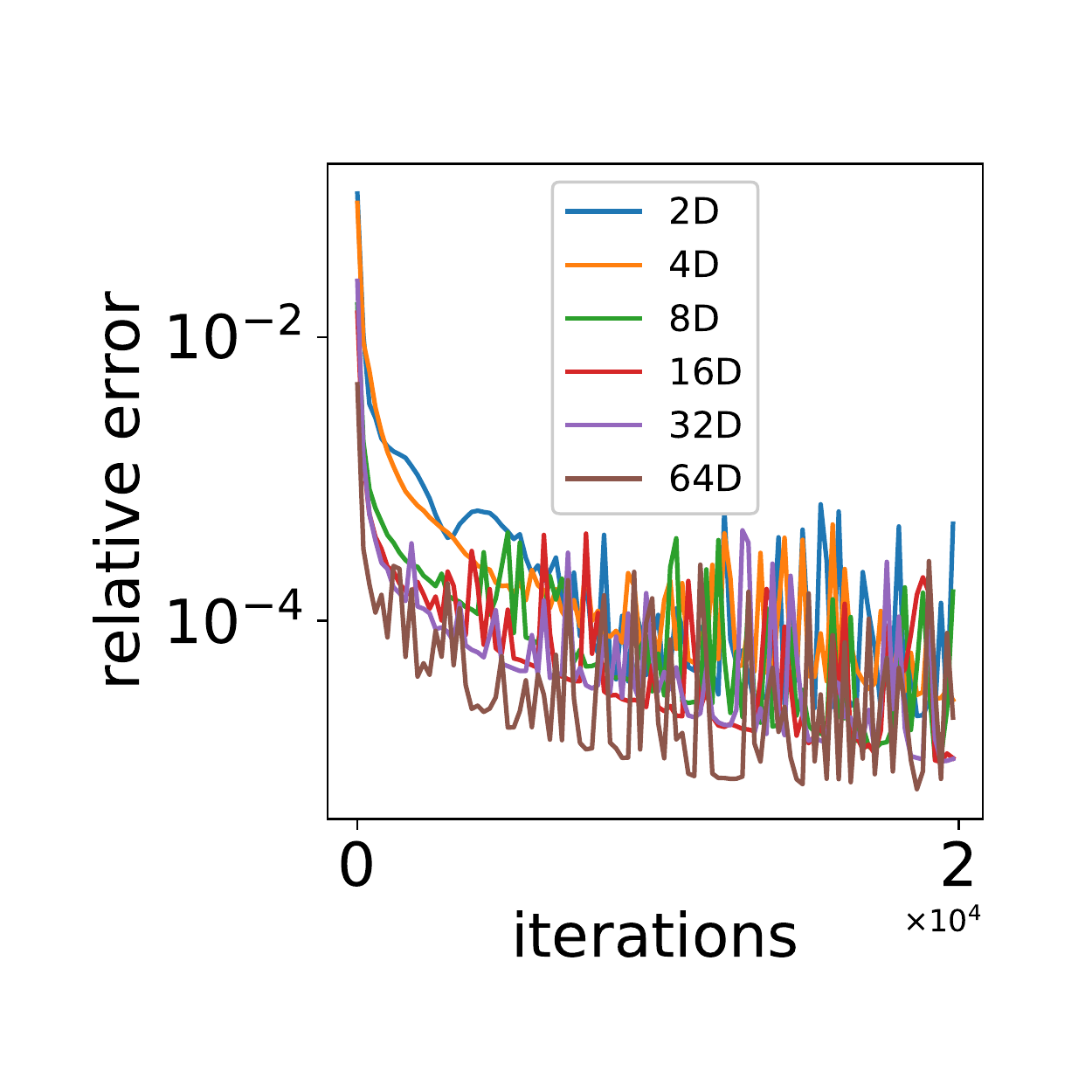}
	}
	\subfigure[MIM \eqref{eqn:rbc new4}-\eqref{eqn:rbc new5}]{
		\includegraphics[width=0.45\textwidth]{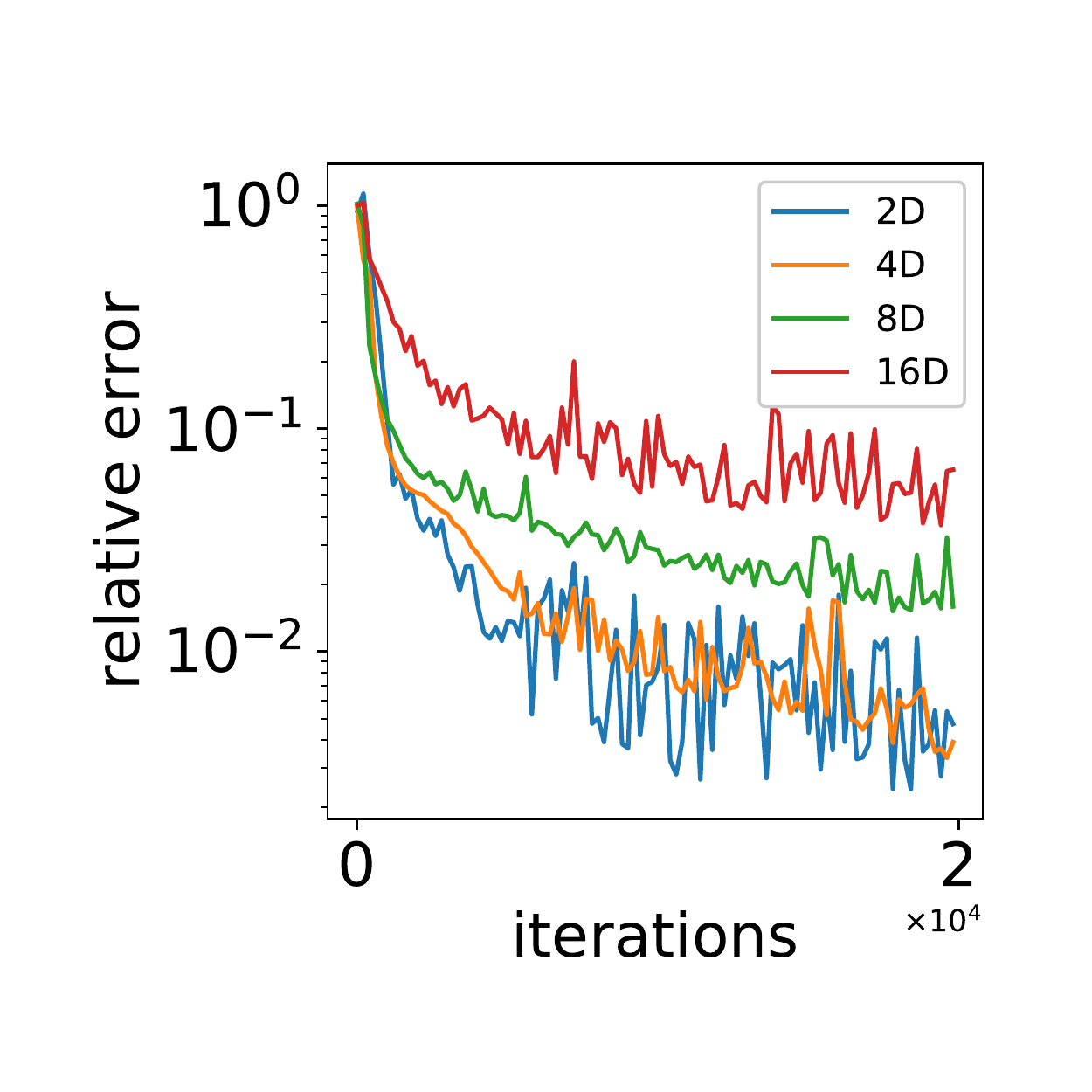}
	}
	\caption{Training processes for Robin boundary condition solved by two formulations of MIM in different dimensions.}
	\label{fig:Robin}
\end{figure}

	The first idea is to introduce $r_1$ and $r_2$ as auxiliary variables to approximate $u\oplus\nabla u , u \ominus \nabla u $ and
	\begin{equation}\label{eqn:rbc new1}
		u = \frac{1}{d}\sum_{i=1}^d\frac{r_{1,i}+r_{2,i}}{2}, \quad p = \frac{r_1 - r_2}{2}.
	\end{equation}
	Then, $r^1_\theta$ and $r^2_\theta$ are constructed to approximate $r_1$ and $r_2$ and satisfy the boundary condition automatically 
	\begin{equation}\label{eqn:rbc new2}
		\begin{aligned}
		&r^1_\theta = x\odot N_\theta(x) \oplus G(x),\\
		&r^2_\theta = (1-x) \odot N^*_\theta(x) \oplus G(x).
		\end{aligned}
	\end{equation}
	Note that both $N_\theta(x)$ and $N^*_\theta(x)$ are $d$ dimensional functions. The corresponding loss function is
	\begin{equation}\label{eqn:rbc new3}
	\begin{aligned}
	L(r_1,r_2) & = \left\|\nabla\left( \frac{1}{d}\sum_{i=1}^d\frac{r_{1,i}+r_{2,i}}{2}\right) -  \frac{r_1 - r_2}{2}\right\|^2_{2,\Omega} \\
	&+ \left\|-\nabla \cdot \frac{r_1 - r_2}{2} +\pi^2 u -f  \right\|^2_{2,\Omega} 
	\end{aligned}
	\end{equation}
	Results of this construction are shown in Table \ref{tbl:Robin} and Table \ref{fig:Robin} with training processes in Figure \ref{fig:Robin}.
	\begin{table}[ht]
		\caption{Relative $L^2$ errors for \eqref{equ:Robin1} with Robin boundary condition solved by MIM with \eqref{eqn:rbc new1}-\eqref{eqn:rbc new3} in different dimensions. ReQu is used as the activation function here and $50000$ sampling points are used in $\Omega$. The training process ends after $50000$ epochs.}
		\label{tbl:Robin}
		\centering\begin{tabular}{|c|c|c|c|}
			\hline
			$d$ & $n$  & $m$ & $\epsilon$ \\
			\hline
			2  & 5 & 2 & 9.47 e-05\\
			4  & 10& 2 & 7.38 e-05\\
			8  & 20& 2 & 4.79 e-05\\
			16 & 20& 2 & 3.80 e-05\\
			32 & 40& 2 & 4.32 e-05\\
			64 & 80& 2 & 3.39 e-05\\
			\hline
		\end{tabular}
	\end{table}

	The second idea is to keep $u$ and use $r$ to represent $u\oplus \nabla u$. The boundary condition can be satisfied by the following construction
	\begin{equation}\label{eqn:rbc new4}
		\begin{aligned}
		&u_\theta = N_\theta(x),\\
		&r_\theta = x\odot(1-x)\odot N^*_\theta(x) \oplus G(x).
		\end{aligned}
	\end{equation}
	Note that $N_\theta(x)$ is a one dimensional function and $N^*_\theta(x)$ is a $d$ dimensional function. The corresponding loss function is defined as 
	\begin{equation}\label{eqn:rbc new5}
	L(u,r) = \|\nabla u \oplus u - r  \|_{2,\Omega}^2 + \| -\nabla \cdot (r \ominus u) + \pi^2 u  - f\|_{2,\Omega}^2
	\end{equation}
	Results are shown in Table \ref{tbl:Robin2}. Compared to Table \ref{tbl:Robin}, one can easily see that the first construction \eqref{eqn:rbc new1}-\eqref{eqn:rbc new3} outperforms the second construction \eqref{eqn:rbc new4}-\eqref{eqn:rbc new5} by more than two orders of magnitude, though the latter has a simpler construction.
	\begin{table}[ht]
		\caption{Relative $L^2$ errors for \eqref{equ:Robin1} with Robin boundary condition solved by MIM with \eqref{eqn:rbc new4}-\eqref{eqn:rbc new5} in different dimensions. ReQu is used as an activation function here and $50000$ sampling points are used in $\Omega$. The training process ends after $20000$ epochs.}
		\centering
		\label{tbl:Robin2}
		\begin{tabular}{|c|c|c|c|}
			\hline
			$d$ & $n$ & $m$ & $\epsilon$\\
			\hline
			2 & 5 & 2 & 7.42 e-03\\
			4 & 10& 2 & 9.71 e-03\\
			8 & 20& 2 & 1.30 e-02\\
			16 &40& 2 & 2.82 e-02\\
			\hline 
		\end{tabular}
	\end{table}
 
	\subsection{Mixed boundary condition}
\begin{figure}
	\centering
	\subfigure[\eqref{eqn:elliptic mbc} with \eqref{eqn:mbc mim1}-\eqref{eqn:mbc mim2}]{
		\includegraphics[width=0.3\textwidth]{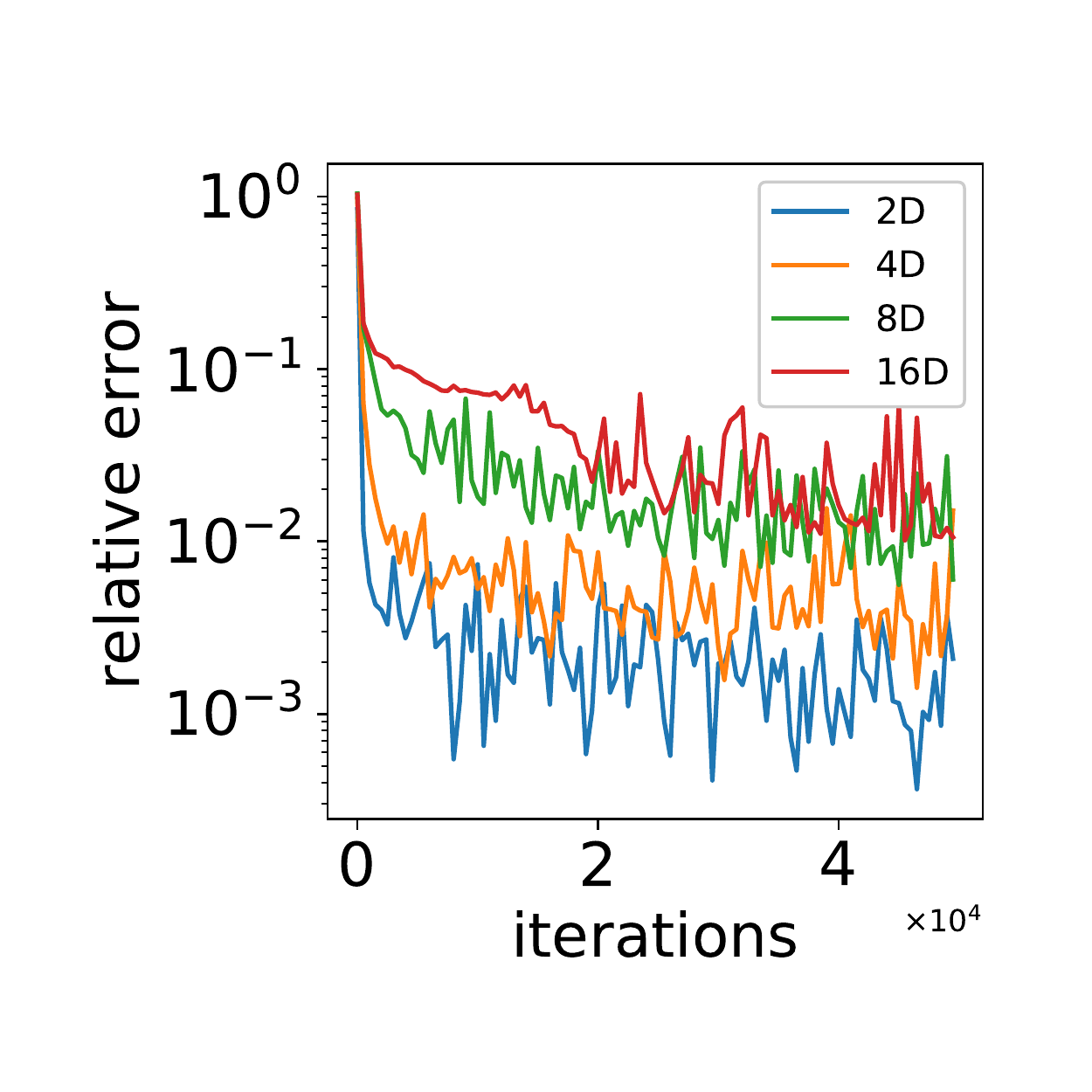}
	}
	\subfigure[\eqref{eqn:elliptic mbc} with \eqref{eqn:mbc mim2}-\eqref{eqn:mbc mim3}]{
		\includegraphics[width=0.3\textwidth]{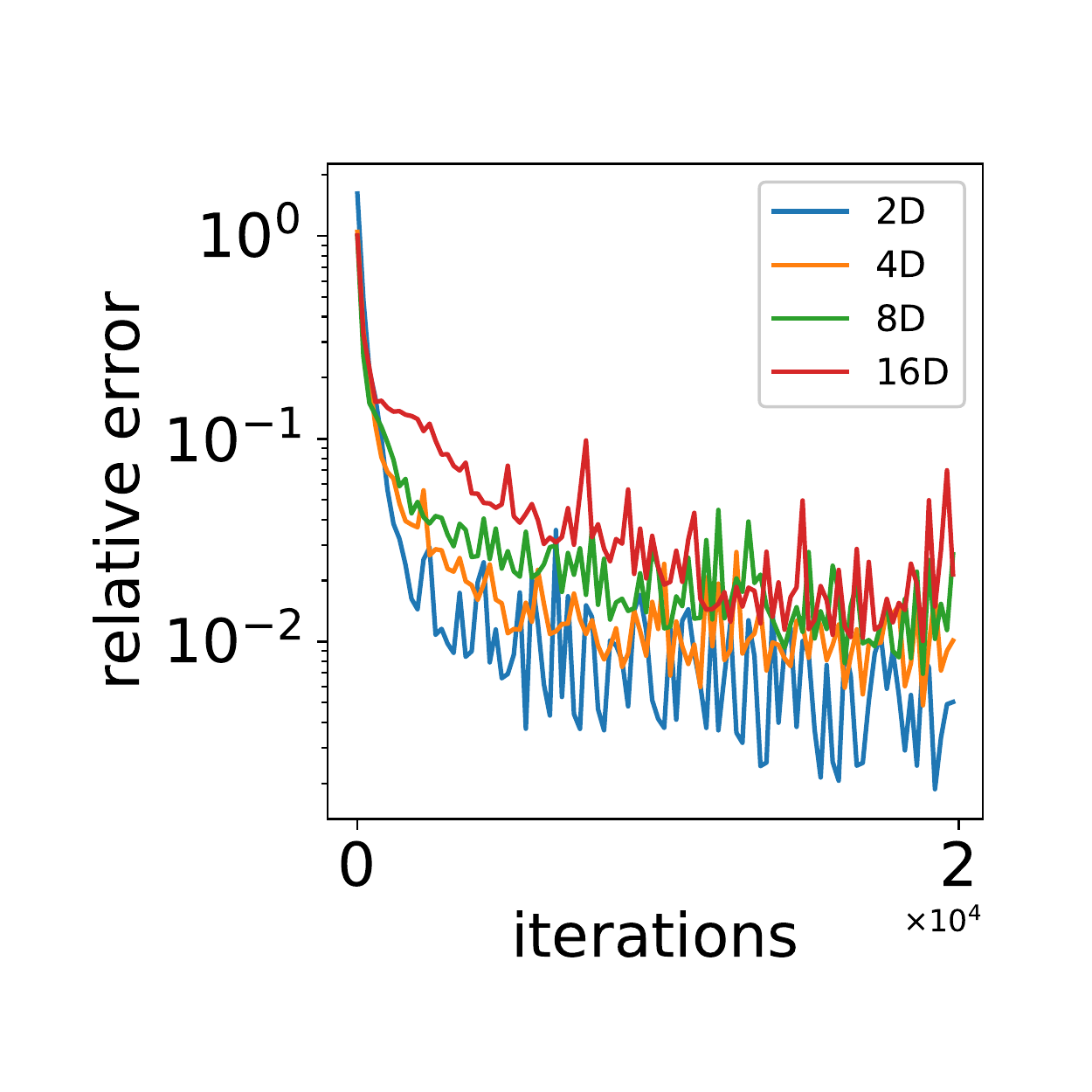}
	}
	\subfigure[\eqref{eqn:elliptic mbc2} with \eqref{eqn:mbc mim2}\&\eqref{eqn:mbc mim4}]{
	\includegraphics[width=0.3\textwidth]{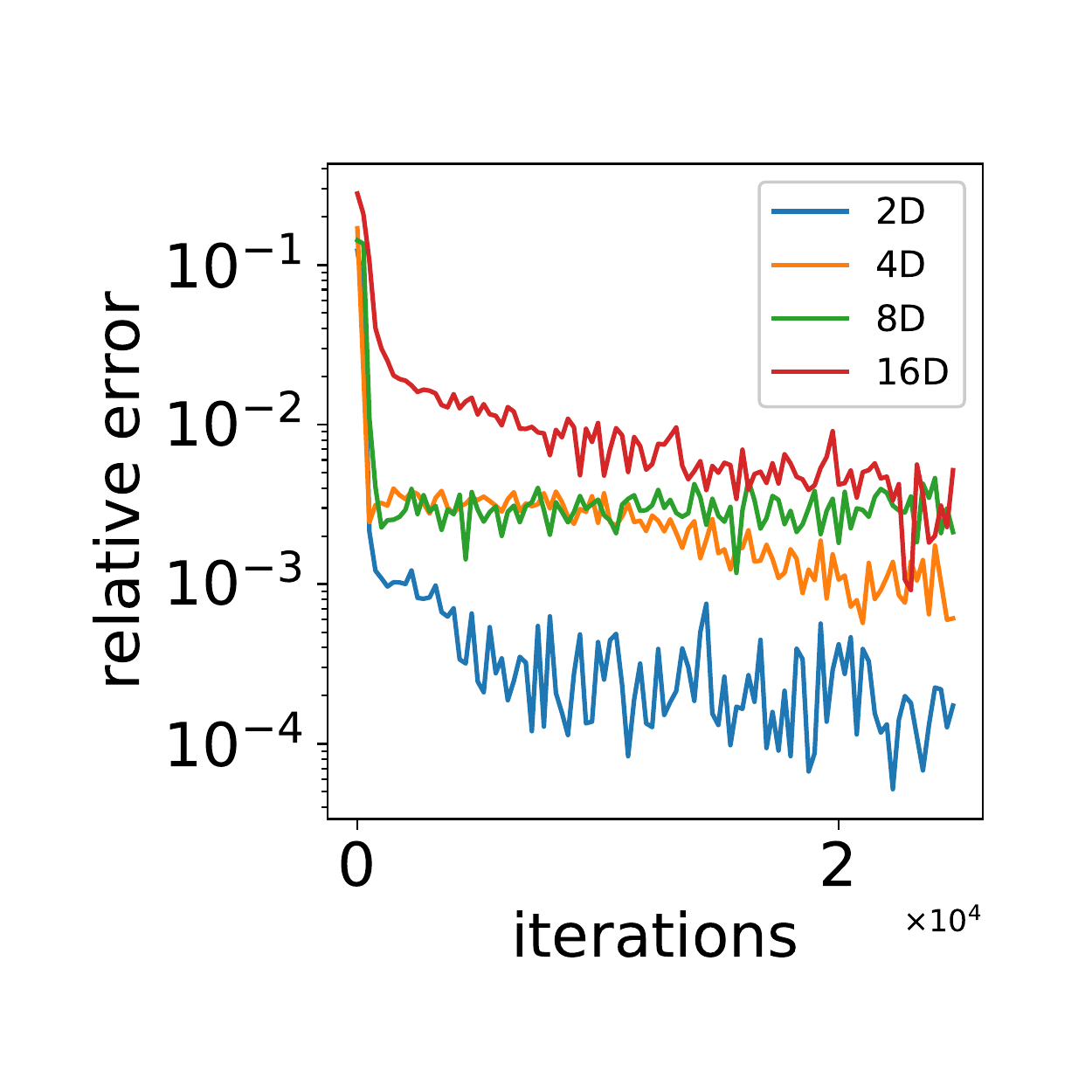}
	}
	\caption{Training processes for \eqref{eqn:elliptic mbc2} with mixed boundary condition solved by MIM  in different dimensions.}
	\label{fig:Mixed}
	\end{figure}

	We consider the Poisson equation with mixed boundary condition
	\begin{equation}\label{eqn:elliptic mbc}
	\left\{
	\begin{aligned}
	&-\Delta u = f 
	& x\in \Omega \\
	&u = 0  	
	& x\in \Gamma_D \\
	&\frac{\partial u}{\partial \nu} = 0 
	&x\in \Gamma_N \\
	\end{aligned}\right.
	\end{equation}
	with exact solution $u(x) = x_1(1-x_1)\sum_{i=2}^d \cos(\pi x_i)$ over $\Omega=[0,1]^d$, $\Gamma_D = \{x: x_1=0\text{ or }x_1 = 1\}$, and $\Gamma_N=\{x: x_i=0 \text{ or } x_i = 1, i \neq 1\}$ . 
	Therefore the trail solution in MIM can be constructed as 
	\begin{equation}\label{eqn:mbc mim1}
		\begin{aligned}
		&\hat{u}_\theta = x_1(1-x_1)N_\theta (x),\\
		& \hat{p}_{\theta,1} = N^*_{\theta,1}(x),\\
		&\hat{p}_{\theta,i} = x_i(1-x_i)N^*_{\theta,i}(x) \quad i \ge 2.
		\end{aligned}
	\end{equation}
	The corresponding loss function is
	\begin{equation}\label{eqn:mbc mim2}
	\begin{aligned}
	L(u,p) = \|\nabla u - p \|_{2,\Omega}^2  + \|\nabla \cdot p + f \|_{2,\Omega}^2.
	\end{aligned}
	\end{equation}
	Results are recorded in Table \ref{tbl:Mixed error} with training processes in Figure \ref{fig:Mixed}(a).
	\begin{table}[ht]
		\centering 
		\caption{Relative $L^2$ errors for \eqref{eqn:elliptic mbc} with mixed boundary condition solved by MIM with construction 
		\eqref{eqn:mbc mim1}-\eqref{eqn:mbc mim2} in different dimensions. ReQu is used as the activation function here and $50000$ sampling points are used in $\Omega$. The training process ends after $50000$ epochs.}
		\begin{tabular}{|c|c|c|c|}
			\hline
			$d$ & $n$ & $m$ &$\epsilon$\\
			\hline
			2  & 5  & 2 & 1.74 e-03\\
			4  & 10 & 2 & 3.87 e-03\\
			8  & 15 & 2 & 1.24 e-02\\
			16 & 24 & 2 & 1.91 e-02\\
			\hline
		\end{tabular}
	\label{tbl:Mixed error}
	\end{table}

	Next we consider \eqref{eqn:elliptic mbc} over a complex domain with exact solution $u(x) =  (x_1-x_2+1) (x_1+x_2)(x_1+\frac{2}{5})x_2(x_2-1)\sum_{i=2}^d \cos(\pi x_i)$, $\Gamma_D = \{x: (x_1-x_2+1) (x_1+x_2)(x_1+\frac{2}{5})x_2(x_2-1)=0\}$, and $\Gamma_N = \{x: x_i=0\text{ or }x_i = 1, i \neq 1,2\}$. $\Omega$ is the interior of the domain formed by $\Gamma_D$ and $\Gamma_N$.
	The trail DNN is constructed as
	\begin{equation}\label{eqn:mbc mim3}
	\begin{aligned}
	&\hat{u}_\theta =  (x_1-x_2+1) (x_1+x_2)(x_1+\frac{2}{5})x_2(x_2-1)N_\theta (x), \\
	& \hat{p}_{\theta,i} = N^*_{\theta,i}(x) \quad  i = 1,2,\\
	&\hat{p}_{\theta,i} = x_i(1-x_i)N^*_{\theta,i}(x) \quad i \ge 3,
	\end{aligned}
	\end{equation}
	and results are shown in Table \ref{tbl:Mixed error2} with training processes in Figure \ref{fig:Mixed}(b).
	\begin{table}[ht]
		\centering 
		\caption{Relative $L^2$ errors for \eqref{eqn:elliptic mbc} with mixed boundary condition solved by MIM with construction 
			\eqref{eqn:mbc mim2}-\eqref{eqn:mbc mim3} in different dimensions. ReQu is used as the activation function here and $50000$ sampling points are used in $\Omega$. The training process ends after $50000$ epochs.}
		\begin{tabular}{|c|c|c|c|}
			\hline
			$d$ & $n$ & $m$ &$\epsilon$\\
			\hline
			2  & 5  & 2 & 5.71 e-03\\
			4  & 10 & 2 & 9.33 e-03\\
			8  & 20 & 2 & 1.35 e-02\\
			16 & 40 & 2 & 1.77 e-02\\
			\hline
		\end{tabular}
		\label{tbl:Mixed error2}
	\end{table}

	Finally we consider the Poisson equation with inhomogeneous mixed boundary condition
	\begin{equation}\label{eqn:elliptic mbc2}
	\left\{
	\begin{aligned}
	&-\Delta u  = f 
	& x\in \Omega \\
	&\frac{\partial u}{\partial \nu} = 0 
	&x \in \Gamma_N = \{x:\|x\|=1\} \\
	& x = g(x) & x \in \Gamma_D = \{x:\|x\|=0.5\} & 
	\end{aligned}\right.
	\end{equation}
	with exact solution $ \cos(\sum_{k=1}^d x_k^2 -1)$ over $\Omega = \{x:0.5<\|x\|<1\}$. As above, the trail solution in MIM is constructed as 
	\begin{equation}\label{eqn:mbc mim4}
	\begin{aligned}
	& \hat{u}_{\theta} = (\sum_{i=1}^d x_i^2 -\frac{1}{4})N_{\theta}(x) + \cos(\frac{3}{4}),\\ 
	& \hat{p}_{\theta} = x\left(- x\cdot N^*_{\theta}(x) \right) + N^*_{\theta}(x),
	\end{aligned}
	\end{equation}
	and results are shown in Table \ref{tab:Laplace Mixed general} with training processes in Figure \ref{fig:Mixed}(c) . By these examples, we see that MIM works well for mixed boundary condition.
	\begin{table}[ht]
		\caption{Relative $L^2$ errors for \eqref{eqn:elliptic mbc2} with mixed boundary condition solved by MIM with construction 
			\eqref{eqn:mbc mim2} and \eqref{eqn:mbc mim4} in different dimensions.}
		\centering
		\begin{tabular}{|c|c|c|c|}
			\hline
			$d$ & $n$ & $m$ & $\epsilon$ \\
			\hline
			2 & 10 & 2 & 2.32 e-04 \\
			4 & 15 & 2 & 8.62 e-04 \\
			8 & 20 & 2 & 2.94 e-03 \\
			16 & 25 & 2 & 3.26 e-03 \\
			\hline
		\end{tabular}
		\label{tab:Laplace Mixed general}
	\end{table}
	\subsection{Periodic boundary condition}
	Consider the elliptic equation over $\Omega = {(-1, 1)}^d$
	\begin{equation}\label{eqn:elliptic pbc}
		-\Delta u + \pi^2 u = f
	\end{equation}
	with periodic boundary condition 
    \begin{equation}\label{eqn:pbc}
			u(x_1,\cdots,x_{i} + I_i,\cdots,x_d) = u(x_1,\cdots,x_{i},\cdots,x_d), \quad i=1,\cdots,d,
	\end{equation}	
    where $I_1 = \cdots = I_d = 2$.
    
	First, we consider the exact solution $u(x) = \sum_{i = 1}^d \left(\cos(\pi x_i)  + \cos(2 \pi x_i)\right) $, and results of MIM are recorded in Table \ref{p0error}.
	Note that the exact solution cannot be explicitly expressed by DNNs when $k=1$. Negligible error is observed if $k\ge 2$ in MIM.
	\begin{table}[H]
		\centering 
		\caption{Relative $L^2$  errors for \eqref{eqn:elliptic pbc} with periodic boundary condition solved by MIM with $k = 1$ in different dimensions. $swish(x)$ is used as the activation function here and $1000$ sampling points are used in $\Omega$. The training process ends after $20000$ epochs when $d=2, 4, 8$ and $50000$ epochs when $d=16$.}
		\begin{tabular}{|c|c|c|c|}
			\hline
			$d$ & $n$ & $m$ & $\epsilon$\\
			\hline
			2  & 8  & 3 & 1.514e-03 \\
			4  & 16 & 3 & 6.593e-03 \\
			8  & 24 & 3 & 1.608e-02 \\
			16 & 32 & 3 & 1.658e-02 \\
			\hline
		\end{tabular}
		\label{p0error}
	\end{table}
	
	Second, we consider the exact solution $u(x) = \sum_{i = 1}^d \cos(\pi x_i) \cos(2 \pi x_i) $, and errors are recorded in Table \ref{p1error}. Note that the exact solution cannot be explicitly expressed by MIM for any $k$.
	\begin{table}[H]
		\centering 
		\caption{Relative $L^2$  errors for \eqref{eqn:elliptic pbc} with periodic boundary condition solved by MIM with $k = 3$ in different dimensions. $swish(x)$ is used as the activation function here and $1000$ sampling points are used in $\Omega$. The training process ends after $20000$ epochs when $d=2, 4, 8$ and $80000$ epochs when $d=16$.}

		\begin{tabular}{|c|c|c|c|}
			\hline
			$d$ & $n$ & $m$ & $\epsilon$\\
			\hline
			2  & 8  & 3 & 2.578e-03 \\
			4  & 8 & 3 & 2.747e-03 \\
			8  & 16 & 3 & 2.965e-03 \\
			16 & 24 & 3 & 3.885e-03 \\
			\hline
		\end{tabular}
		\label{p1error}
	\end{table}
	
\begin{figure}[H]
	\centering
	\includegraphics[scale = 0.8]{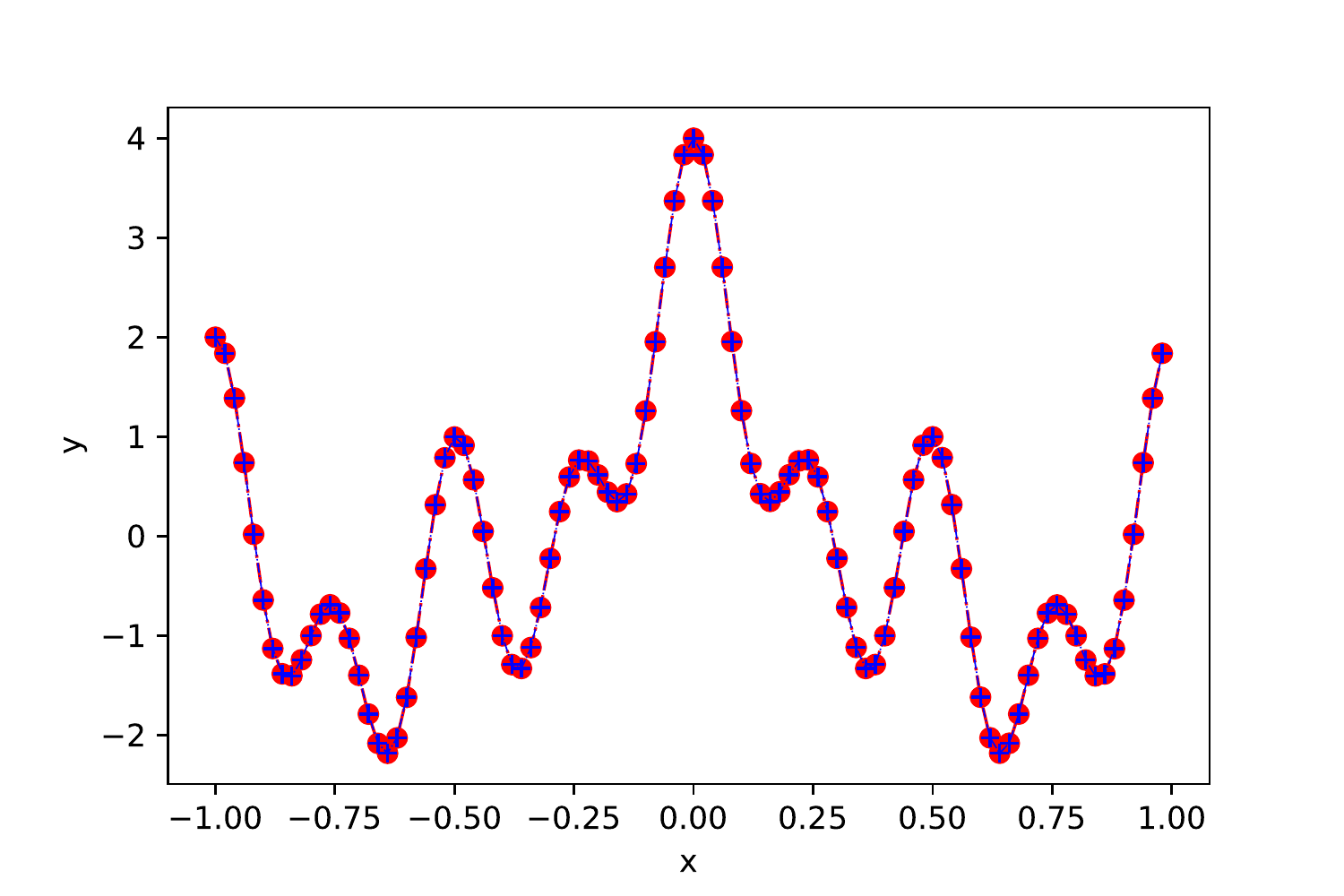}
	\caption{Exact (line) and neural (circle) solutions with high-frequency information in 1D. Parameters used here are $k = 1$, $n = 20$, and $m = 3$. $swish(x)$ is used as the activation function here and $1000$ sampling points are used in $\Omega$. The training process ends after $20000$ epochs.}
\label{fig: 1dsolution}
	\end{figure}

	Finally, we consider a 1D solution $y =  \cos(\pi x) + \cos(2 \pi x) + \cos(4 \pi x) + \cos(8 \pi x)$ with high-frequency information.
	The relative $L^2$ error after training is $0.0043$. The exact and neural solutions are plotted in Figure \ref{fig: 1dsolution}. We shall not expect a good approximation at the first glance since the high-frequency component cannot be expressed when $k=1$. However, due to the deep nature of neural networks used in MIM, the high-frequency component is actually resolved well by MIM; see Figure \ref{fig: 1dsolution}. We also observe that a larger $k$ is needed in high dimensions. 

	\section{Enforcement of exact initial conditions and numerical results}
	\label{sec:numerical result initial}
	
   	\subsection{Parabolic equation}
   
	Consider the parabolic equation with initial and boundary conditions
	\begin{equation}\label{eqn:parabolic}
	\left\{
	\begin{aligned}
	& u_{t} - \Delta u = f & (t, x) \in  (0, 1)\times \Omega\\
	& u(x,t) = 0 &  (t, x) \in  (0, 1)\times\partial \Omega \\
	& u(x, 0) = 0 & x \in  \Omega 
	\end{aligned}\right.
	\end{equation}
    with the exact solution $u(x,t) = t  \Pi_{i=1}^d \sin(\pi x_i)$ over $ \Omega  = {(0, 1)}^d$.
    
\begin{figure}
	\centering
	\subfigure[2D]{
		\includegraphics[width=0.3\textwidth]{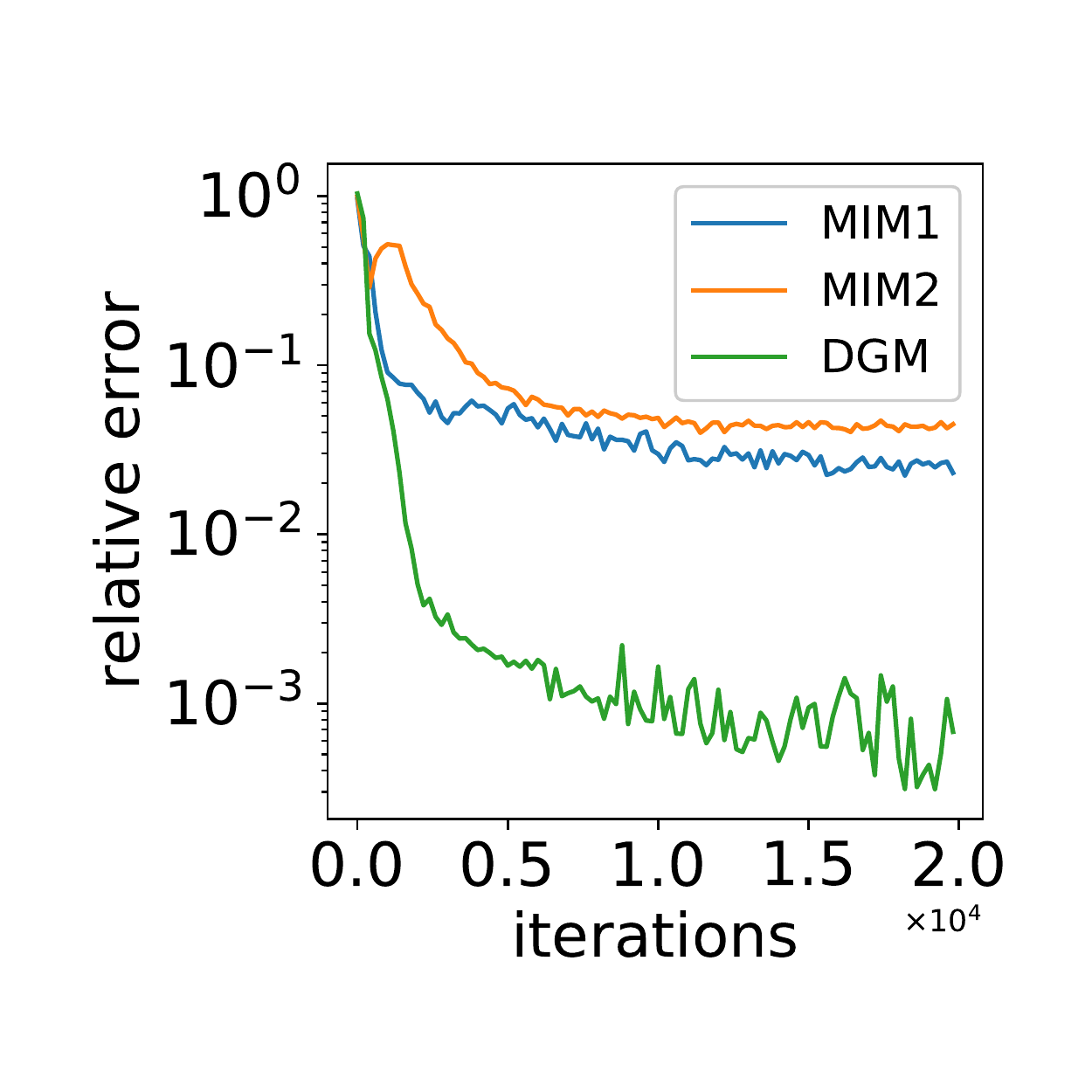}
	}
	\subfigure[3D]{
		\includegraphics[width=0.3\textwidth]{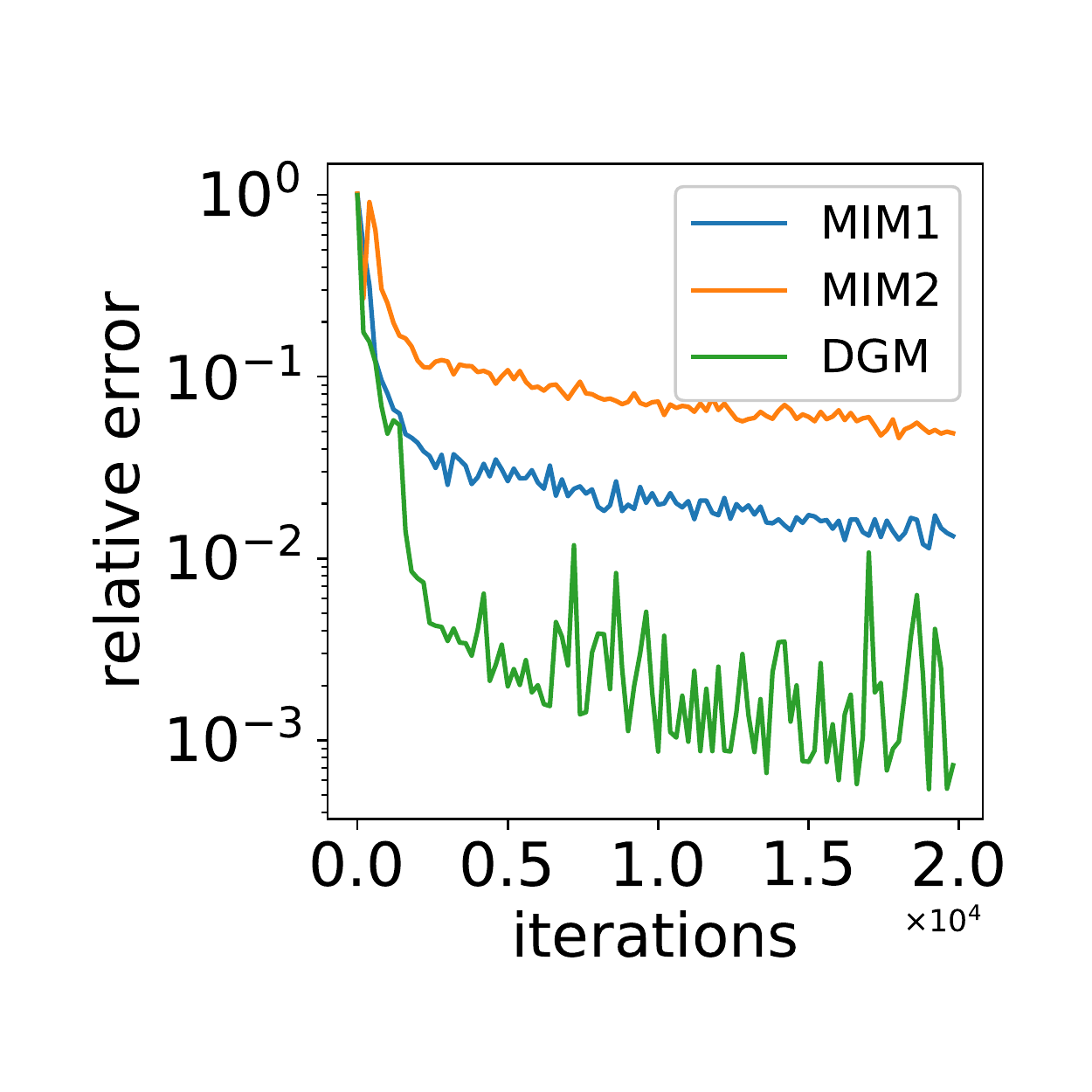}
	}
	\subfigure[5D]{
		\includegraphics[width=0.3\textwidth]{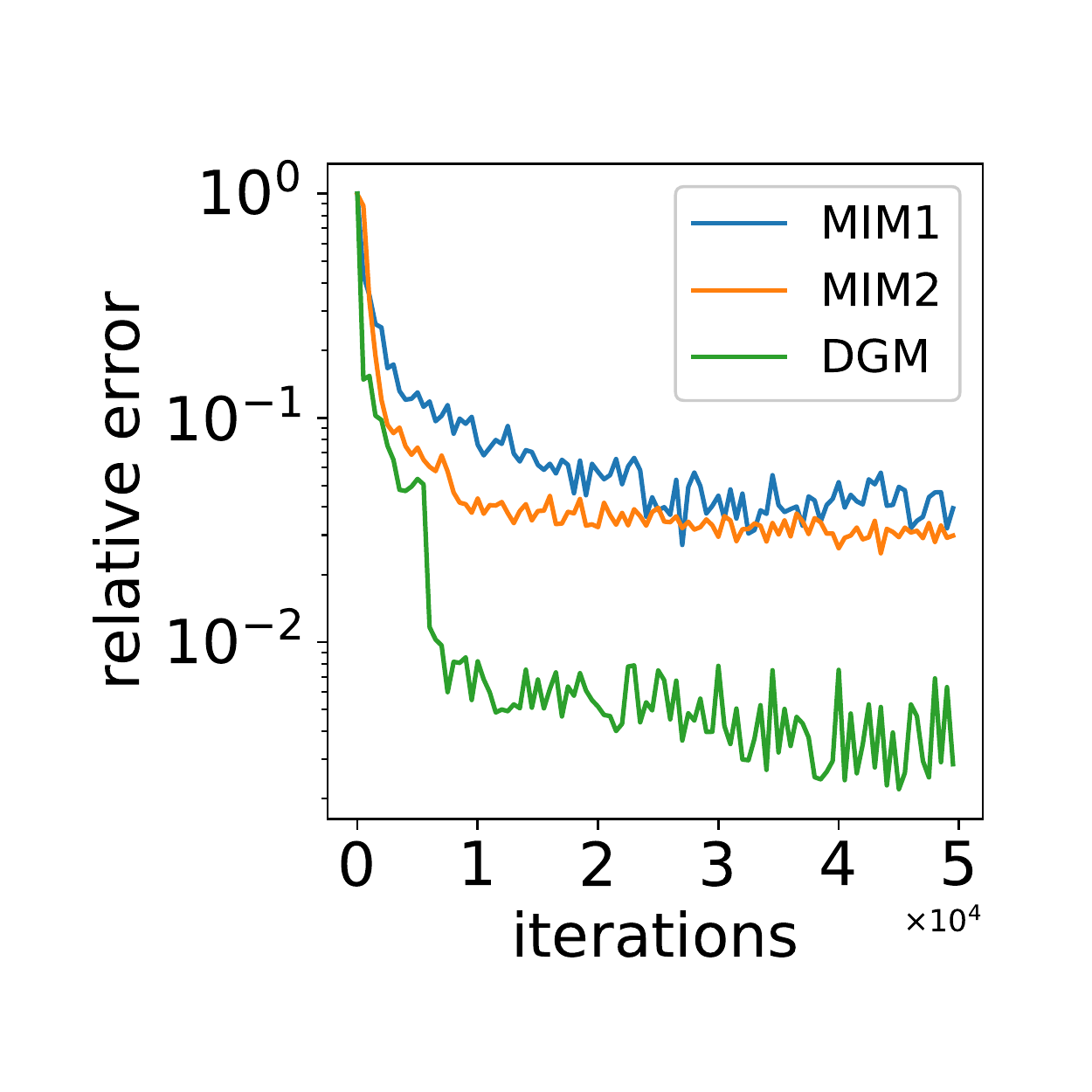}
	}
	\subfigure[10D]{
		\includegraphics[width=0.3\textwidth]{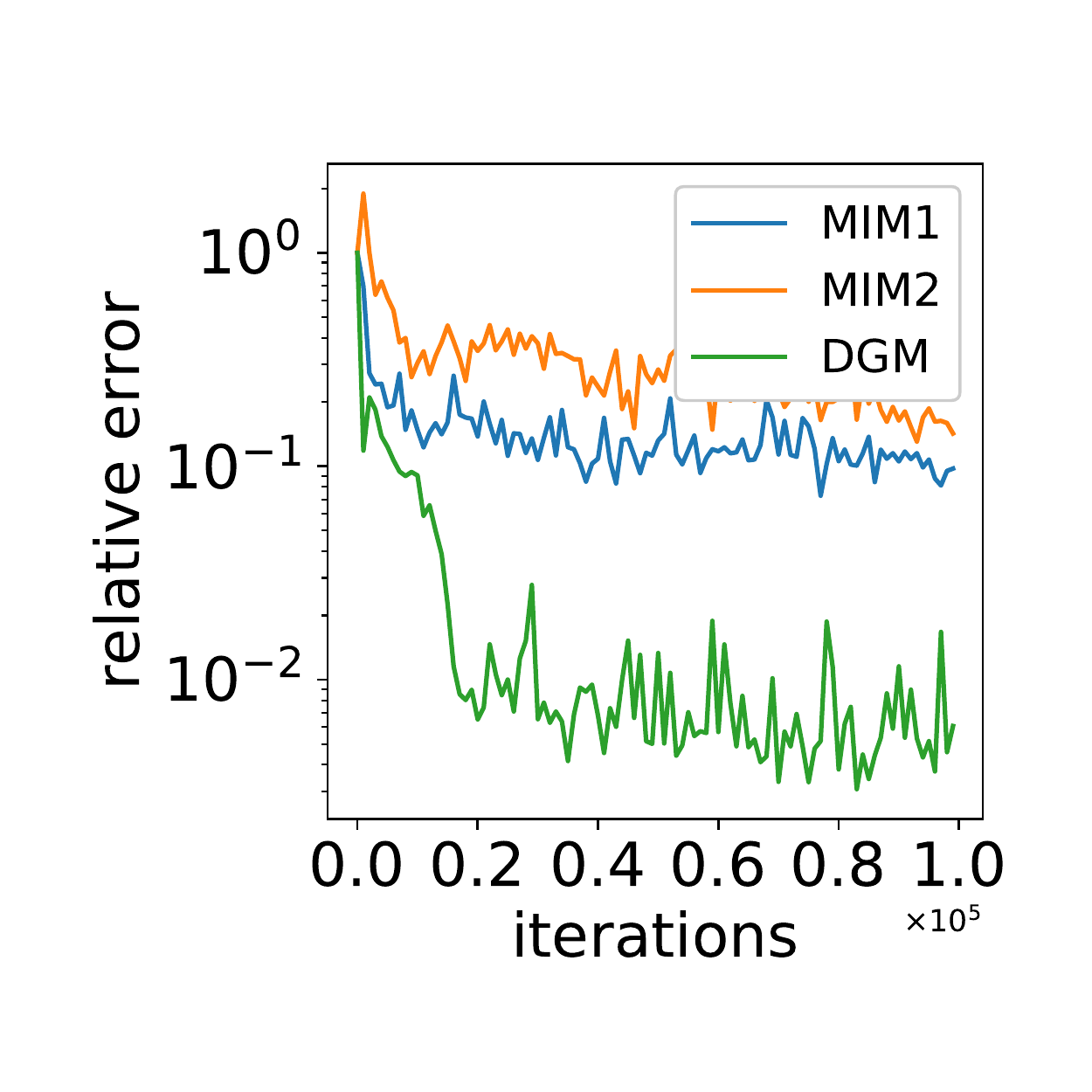}
	}
	\subfigure[12D]{
		\includegraphics[width=0.3\textwidth]{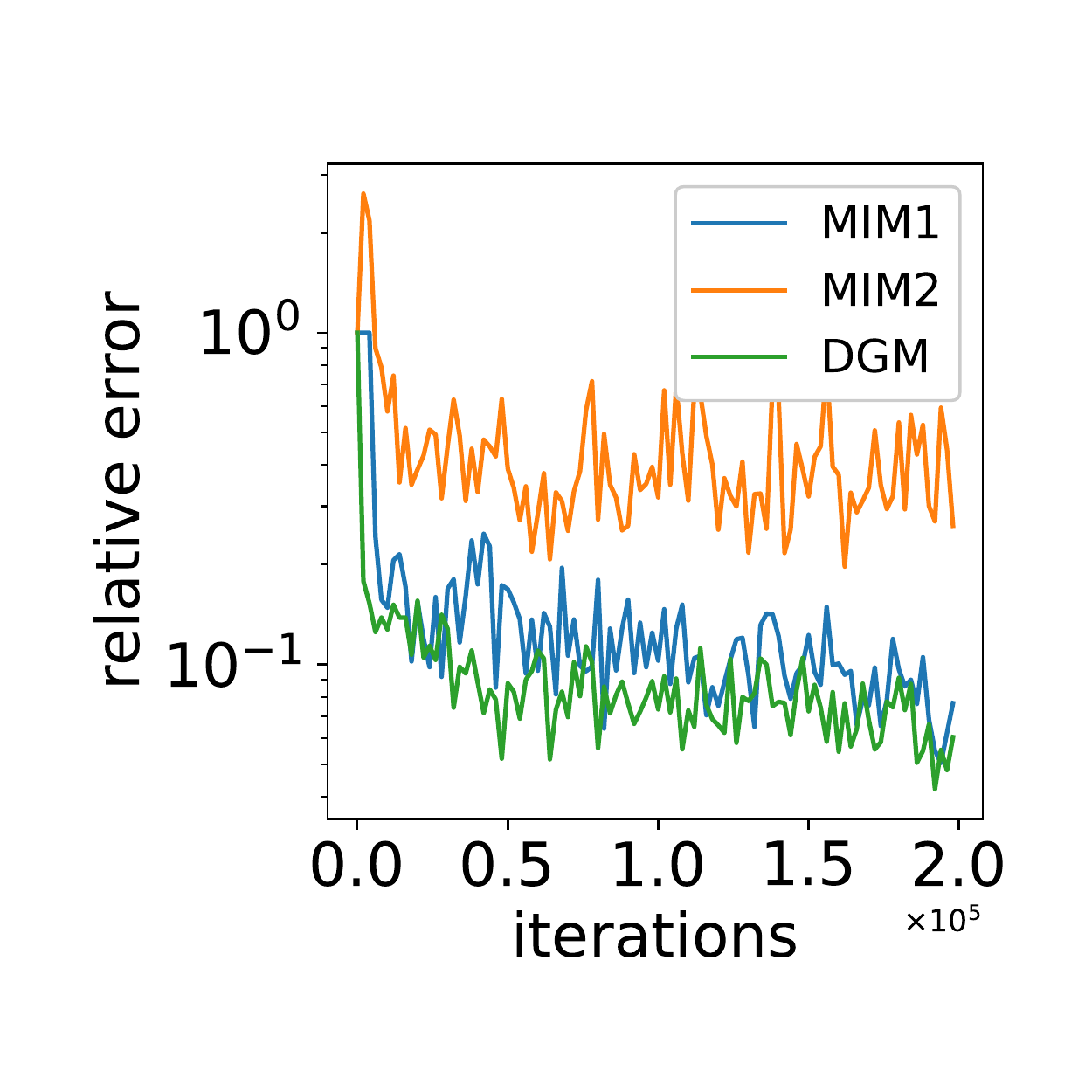}
	}
	\caption{Training processes for parabolic equation \eqref{eqn:parabolic} by DGM, MIM1 and MIM2 in different dimensions.}
	\label{fig:Parabolic}
	\end{figure}

    Both initial and boundary conditions are imposed directly on the solution, therefore we can easily construct trail DNNs that satisfy exact conditions for both DGM and MIM in the following form
    \begin{equation}\label{eqn:parabolic ic1}
    	u_\theta(x, t) = t\Pi_{i=1}^d(x_i-x_i^2) N_{\theta}(x, t). 
    \end{equation}
    The corresponding loss function in DGM is 
    \begin{equation}\label{eqn:parabolic ic2}
	L(u) = \|u_{t} -\Delta u -f \|^2_{L^2((0,1)\times\Omega)} 
	\end{equation}
    There are two options for the loss function in MIM for \eqref{eqn:parabolic}: MIM1 and MIM2. They are
    \begin{equation}\label{eqn:parabolic ic mim1}
	\begin{aligned}
	L(u,v,p) & =  \|v-\nabla \cdot p -f \|^2_{L^2((0,1)\times\Omega)} + \|p-\nabla u\|^2_{L^2((0,1)\times\Omega)} \\
	&+ \|v - u_{t}\|^2_{L^2((0,1)\times\Omega)}
	\end{aligned}    
    \end{equation}
   and
    \begin{equation}\label{eqn:parabolic ic mim2}
		L(u,v,p) =  \|u_{t}-\nabla \cdot p -f \|^2_{L^2((0,1)\times\Omega)} + \|p-\nabla u\|^2_{L^2((0,1)\times\Omega)},
    \end{equation}
    respectively. Results of DGM, MIM1, and MIM2 are recorded in Table \ref{paraerror}. Since all three methods are free of penalty terms, they all work well for \eqref{eqn:parabolic}. As the dimension $d$ increases, MIM starts to outperform DGM. Figure \ref{fig:Parabolic} plots the training processes of DGM, MIM1, and MIM2 in different dimensions.
	\begin{table}[ht]
		\centering 
		\caption{Relative $L^2$ errors of DGM, MIM1, and MIM2 for parabolic equation \eqref{eqn:parabolic} in different dimensions. $swish(x)$ is used as the activation function here and $2000$ sampling points are used in $\Omega$. The training process ends after $50000$ epochs when $d=2, 3, 5$, $100000$ epochs when $d=10$, and $200000$ epochs when $d=12$.}
			\begin{tabular}{|c|c|c|c|c|c|}
			\hline
			\multirow{2}*{$d$} & \multirow{2}*{$n$} & \multirow{2}*{$m$}  & \multicolumn{3}{c|}{$\epsilon$} \\
			\cline{4-6}
			~& ~ & ~  & MIM1 & MIM2 & DGM\\
			\hline
			2  & 4 & 3 & 1.92 e-02 & 4.27 e-02 & 5.16 e-04 \\
			3  & 8 & 3 & 1.42 e-02 & 3.83 e-02 &1.74 e-04\\
            5  & 8 & 3 & 3.48 e-02& 3.22 e-02&1.49 e-03\\
			10  & 20 & 3 & 8.17 e-02 & 1.32 e-01 &4.70e-03\\
            12  & 20 & 3 & 7.47 e-02 & 2.20 e-01&5.06e-02\\
			\hline
		\end{tabular}
     	\label{paraerror}
	\end{table}

\subsection{Wave equation}

\begin{figure}
		\centering
		\subfigure[$n = 10$, ReQu]{
			\includegraphics[width=0.3\textwidth]{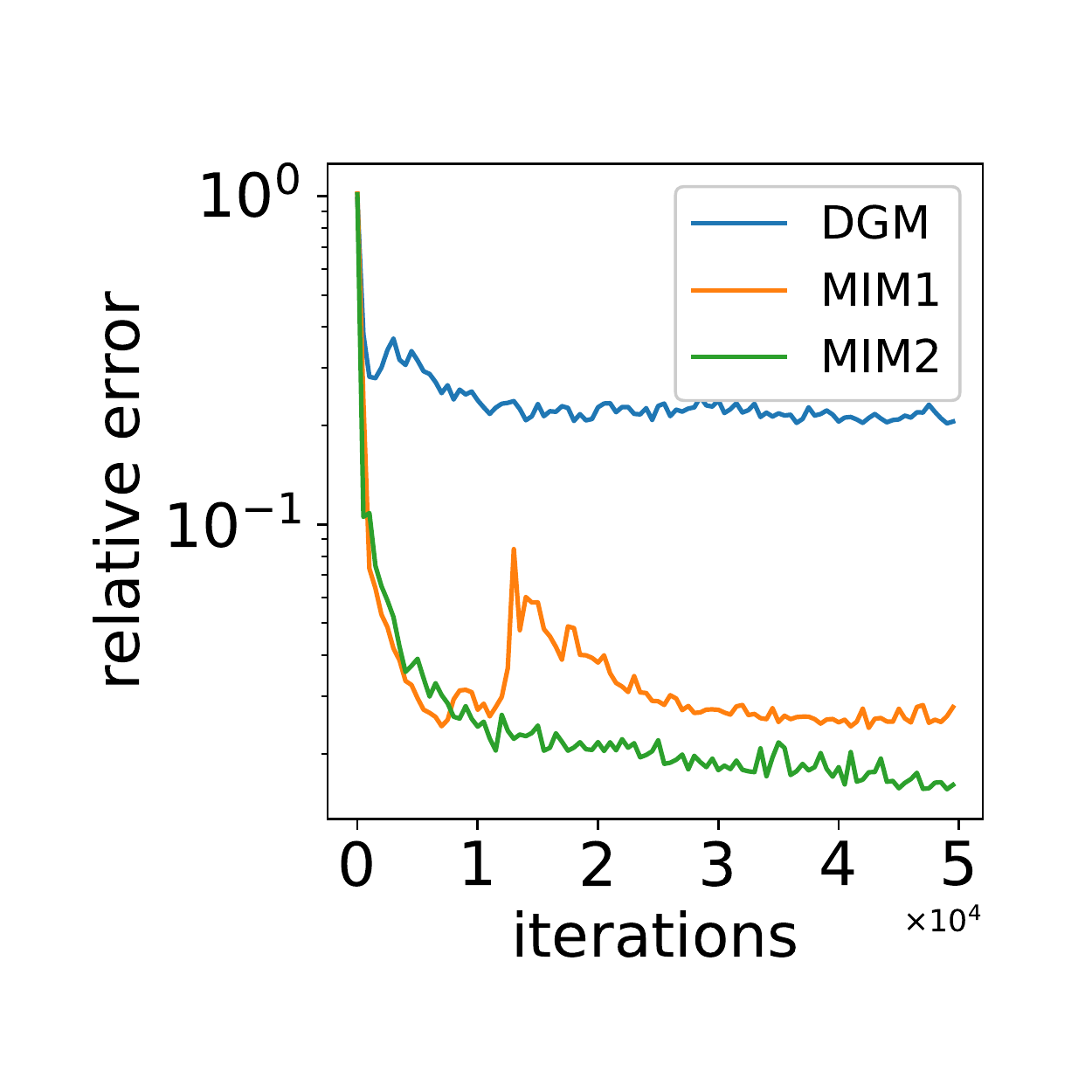}
		}
		\subfigure[$n = 20$, ReQu]{
			\includegraphics[width=0.3\textwidth]{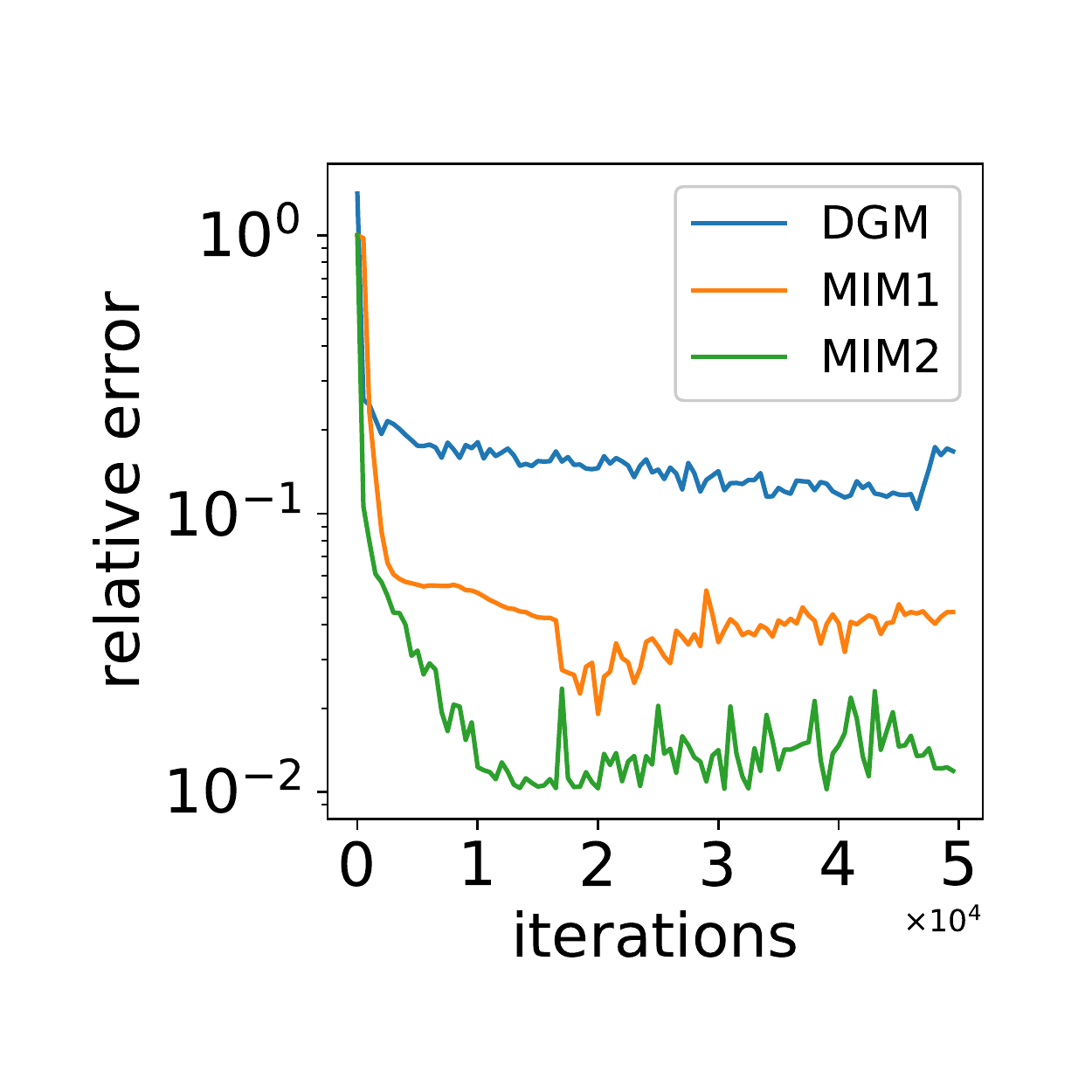}
		}
		\subfigure[$n = 40$, ReQu]{
			\includegraphics[width=0.3\textwidth]{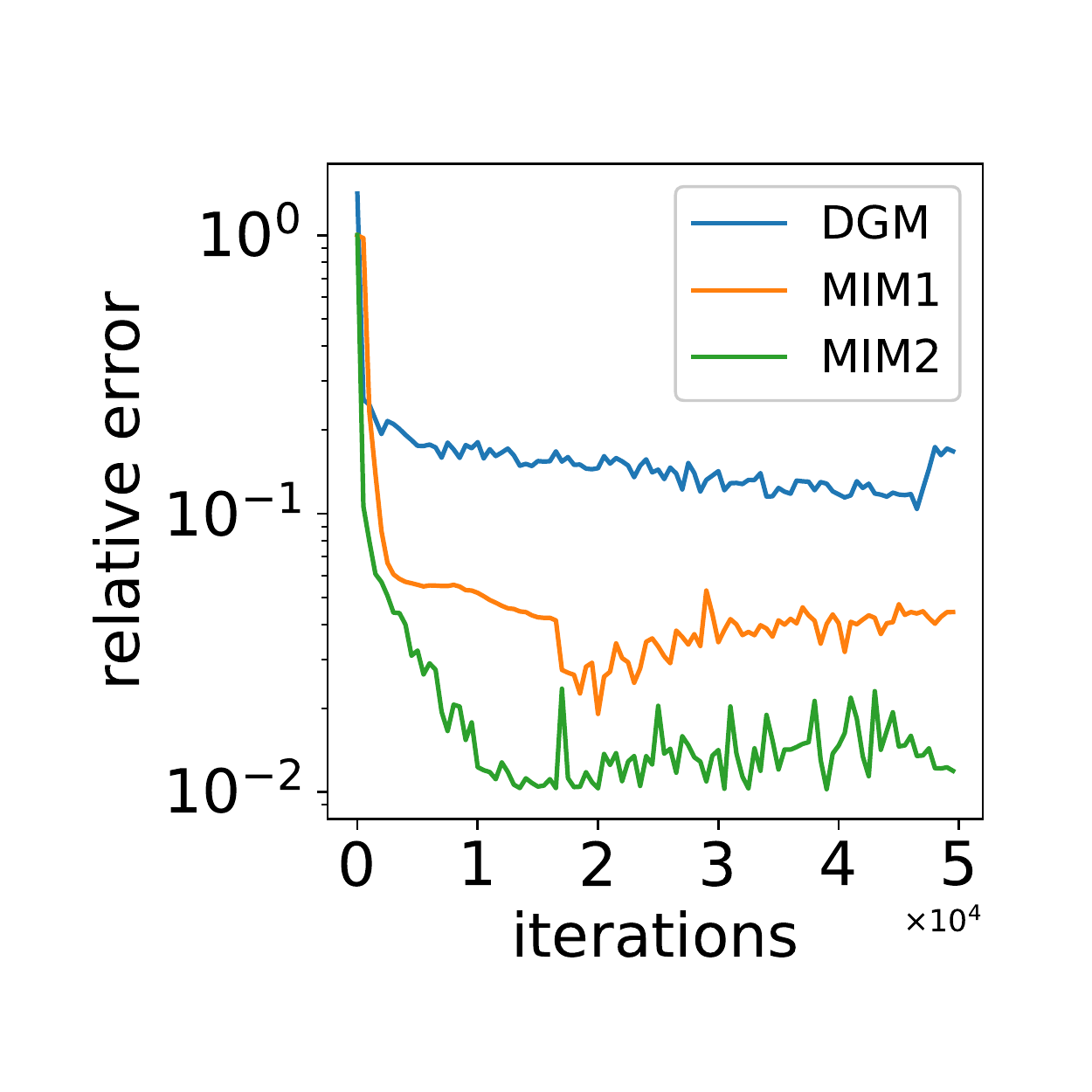}
		}
		\subfigure[$n = 10$, ReCu]{
			\includegraphics[width=0.3\textwidth]{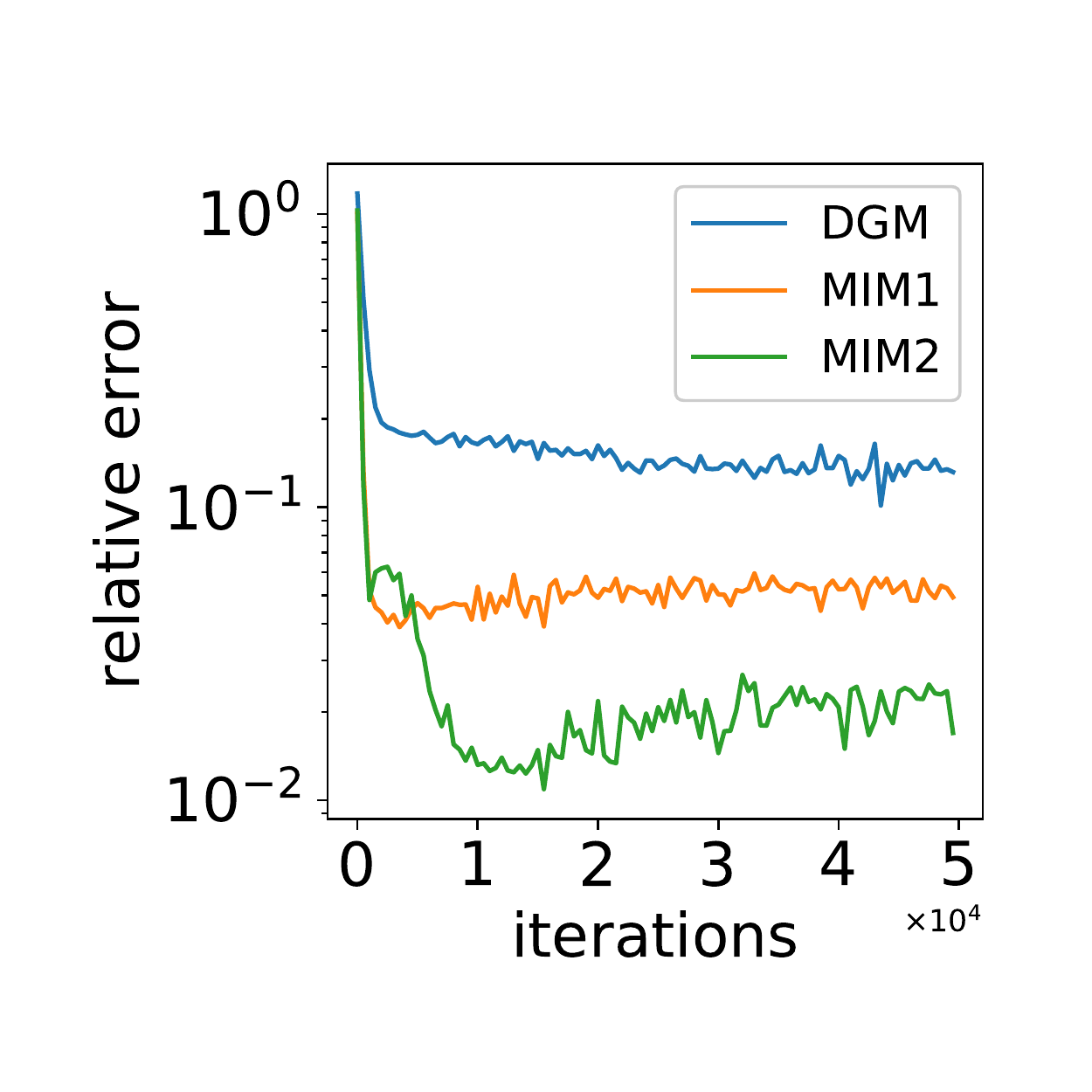}
		}
		\subfigure[$n = 20$, ReCu]{
			\includegraphics[width=0.3\textwidth]{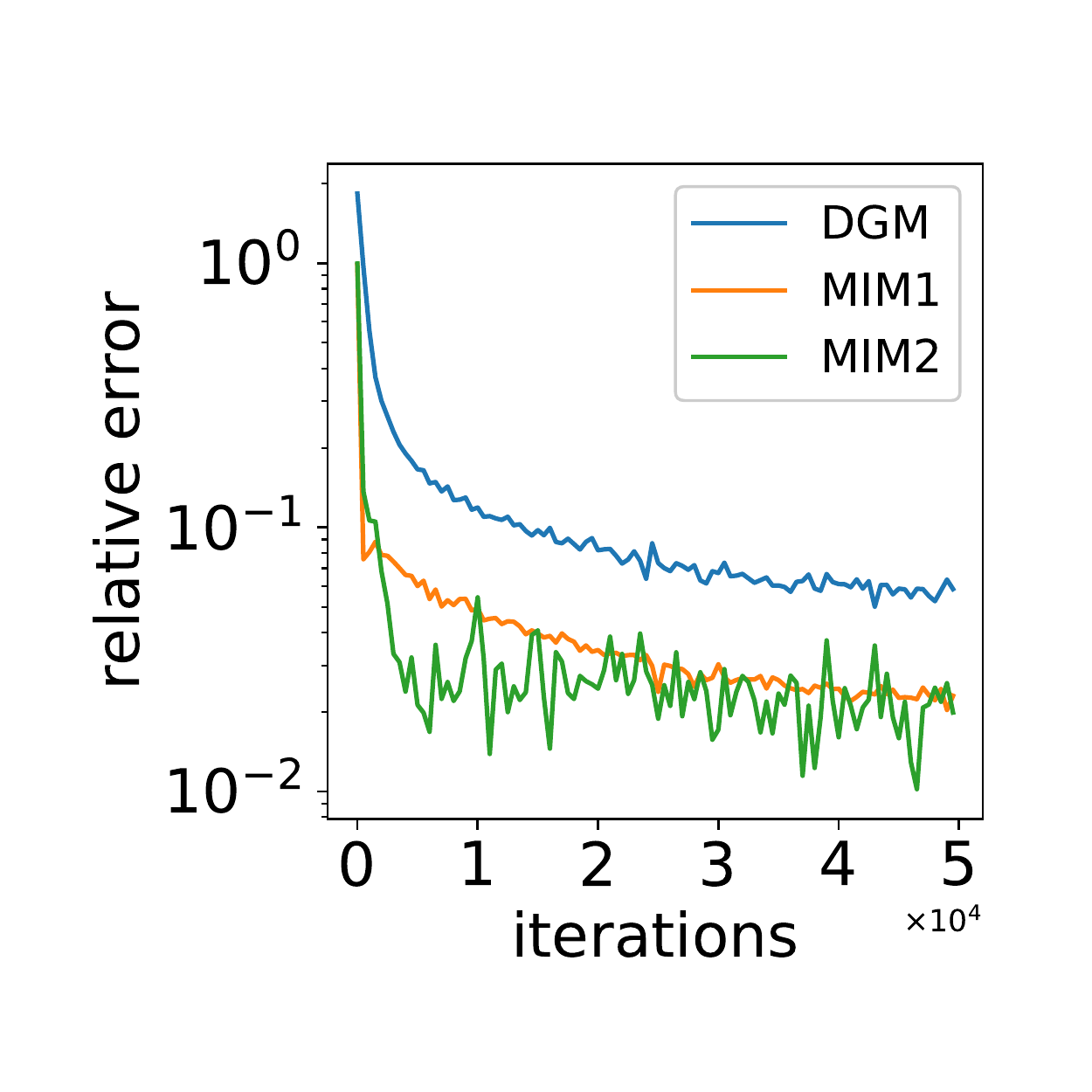}
		}
		\subfigure[$n = 40$, ReCu]{
			\includegraphics[width=0.3\textwidth]{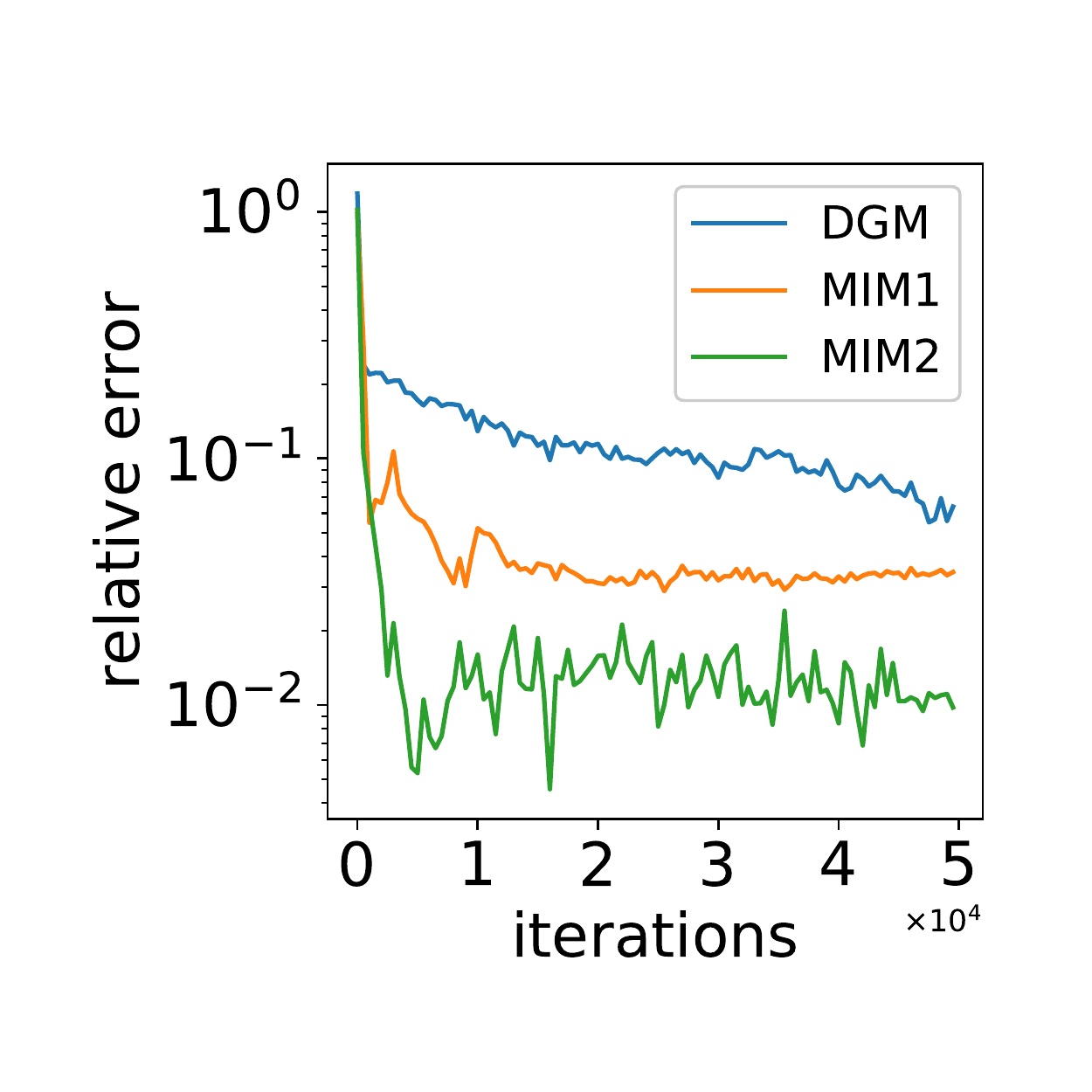}
		}
		\caption{Training processes for wave equation \eqref{eqn:wave} by DGM, MIM1 and MIM2 with the network depth $m = 3$ when $d=2$.}
		\label{fig:wave 3D2L}
\end{figure}

\begin{figure}
	\centering
	\subfigure[$n = 10$, ReQu]{
		\includegraphics[width=0.3\textwidth]{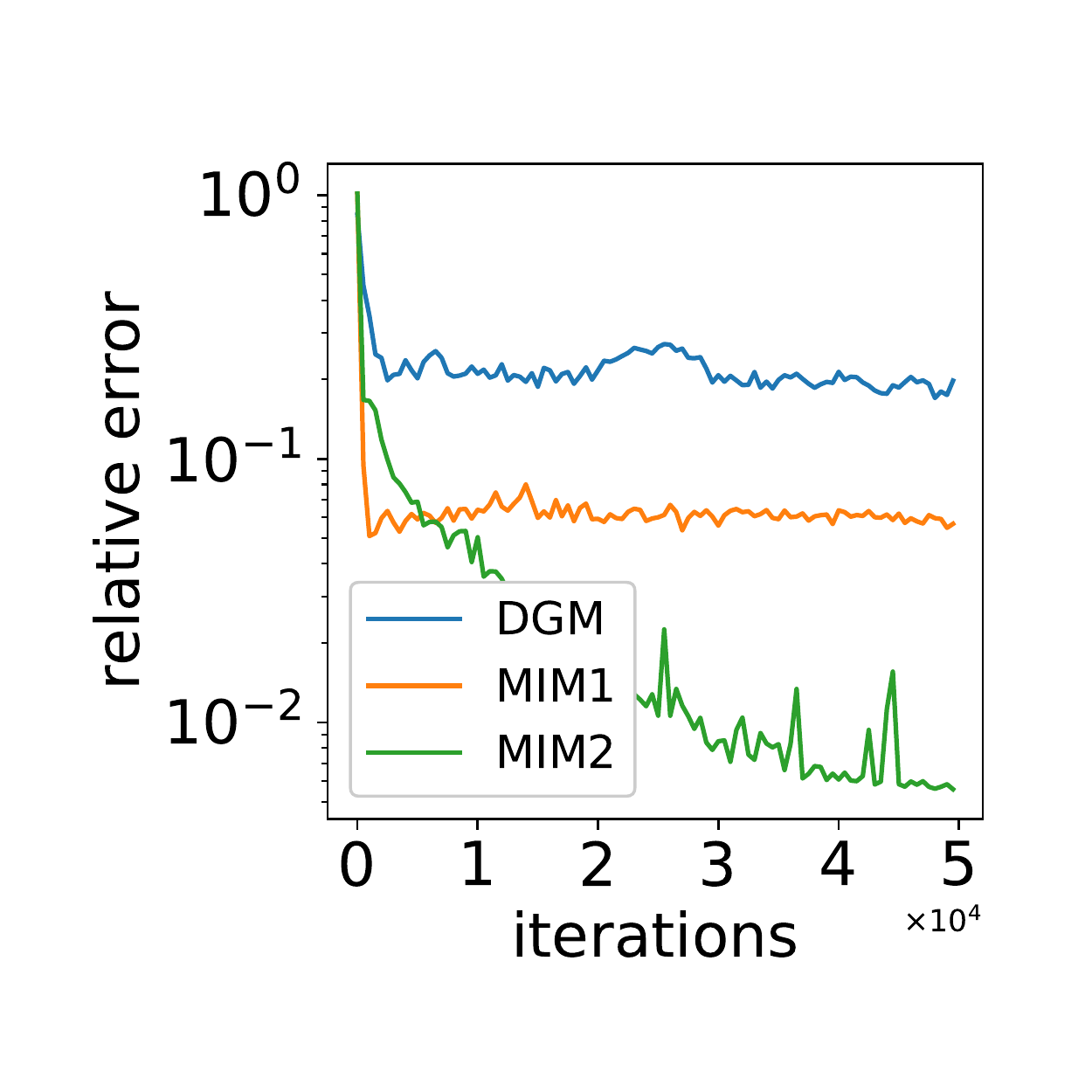}
	}
	\subfigure[$n = 20$, ReQu]{
		\includegraphics[width=0.3\textwidth]{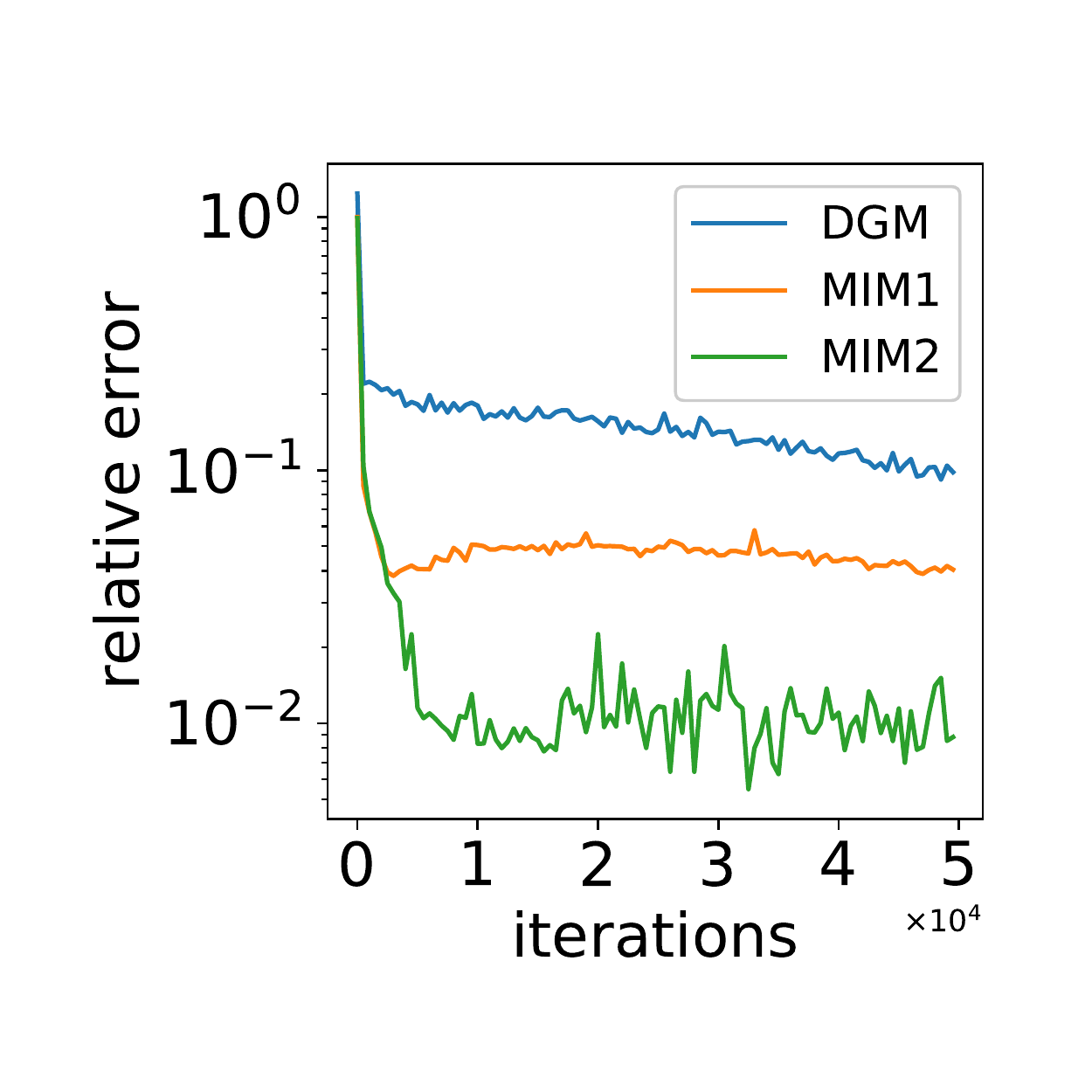}
	}
	\subfigure[$n = 40$, ReQu]{
		\includegraphics[width=0.3\textwidth]{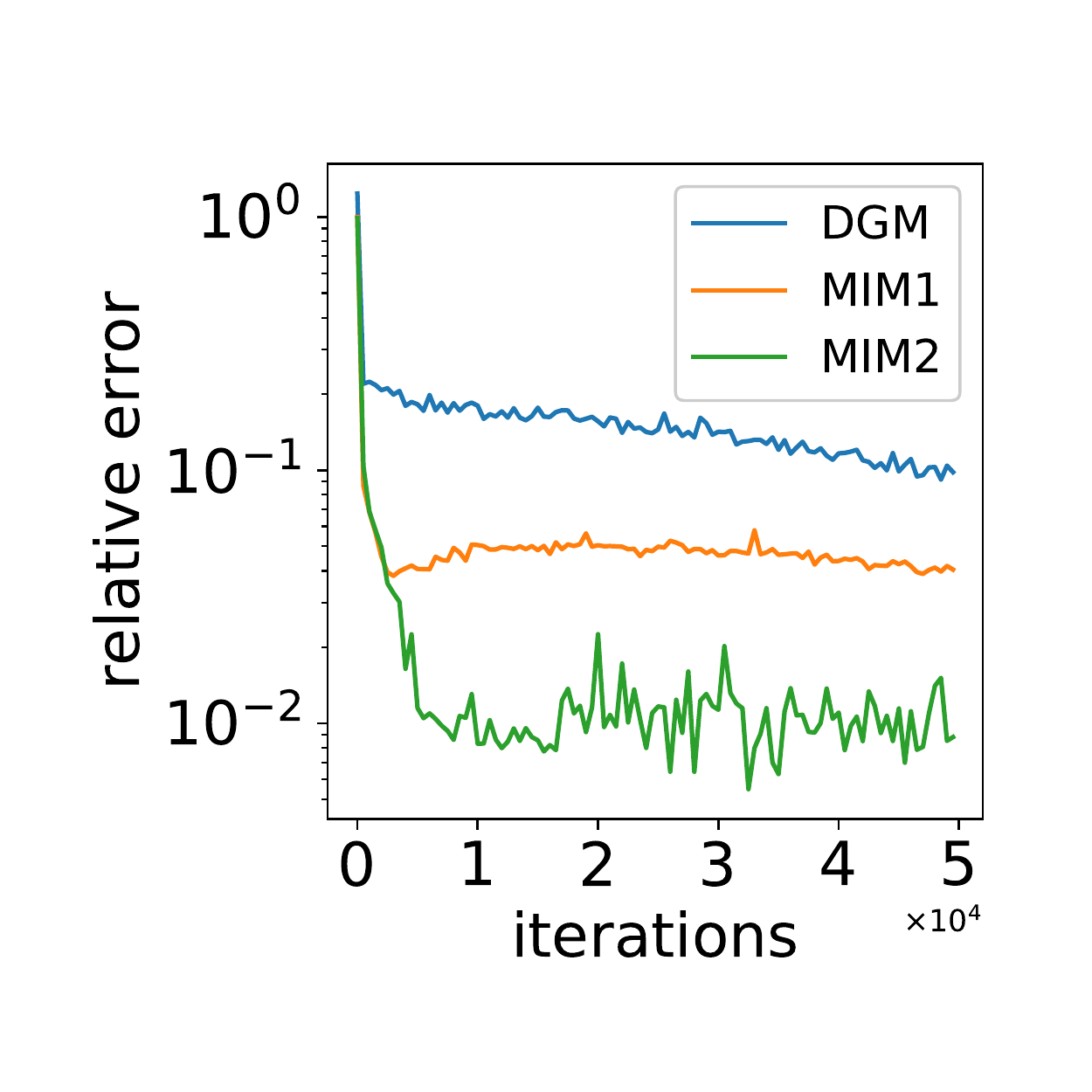}
	}
	\subfigure[$n = 10$, ReCu]{
		\includegraphics[width=0.3\textwidth]{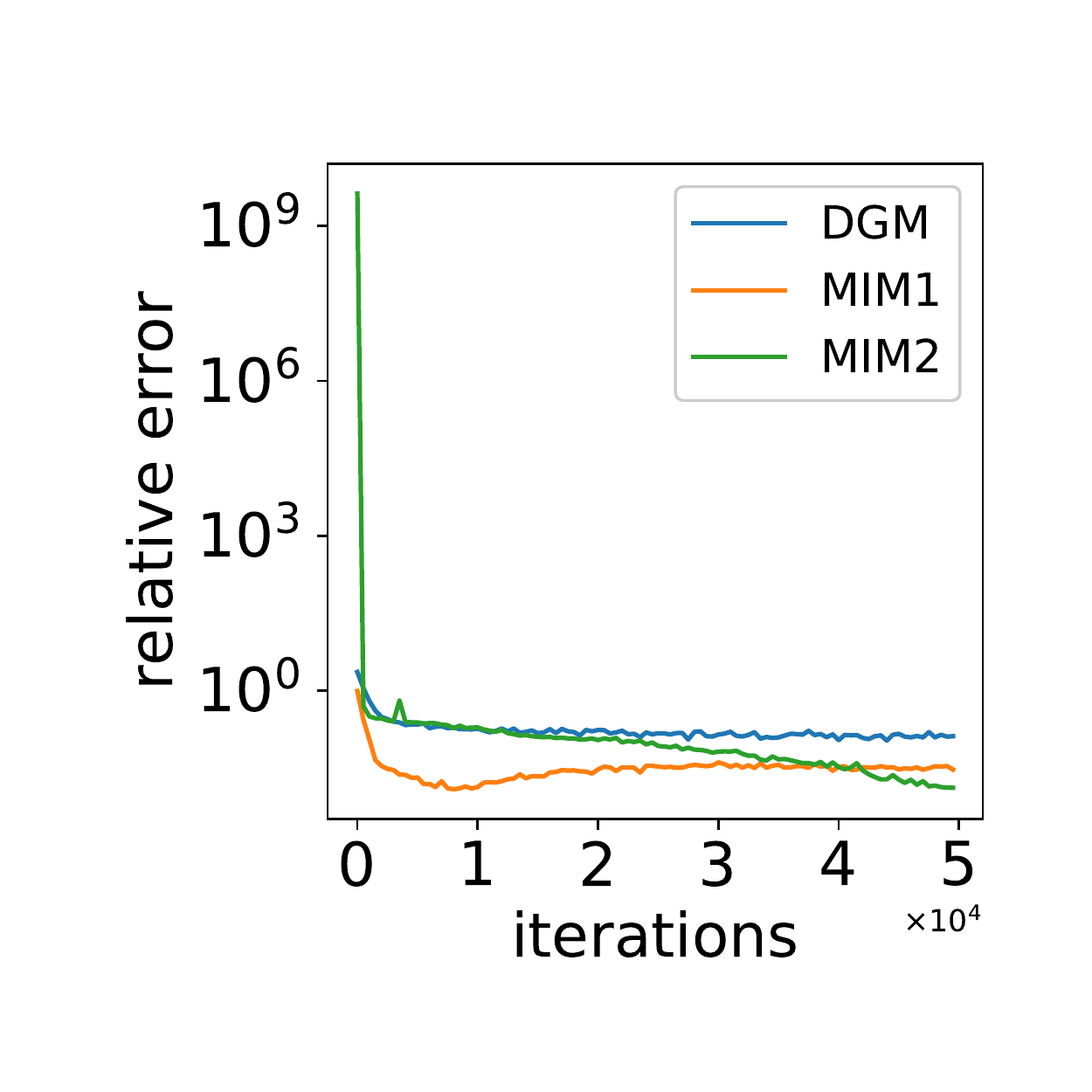}
	}
	\subfigure[$n = 20$, ReCu]{
		\includegraphics[width=0.3\textwidth]{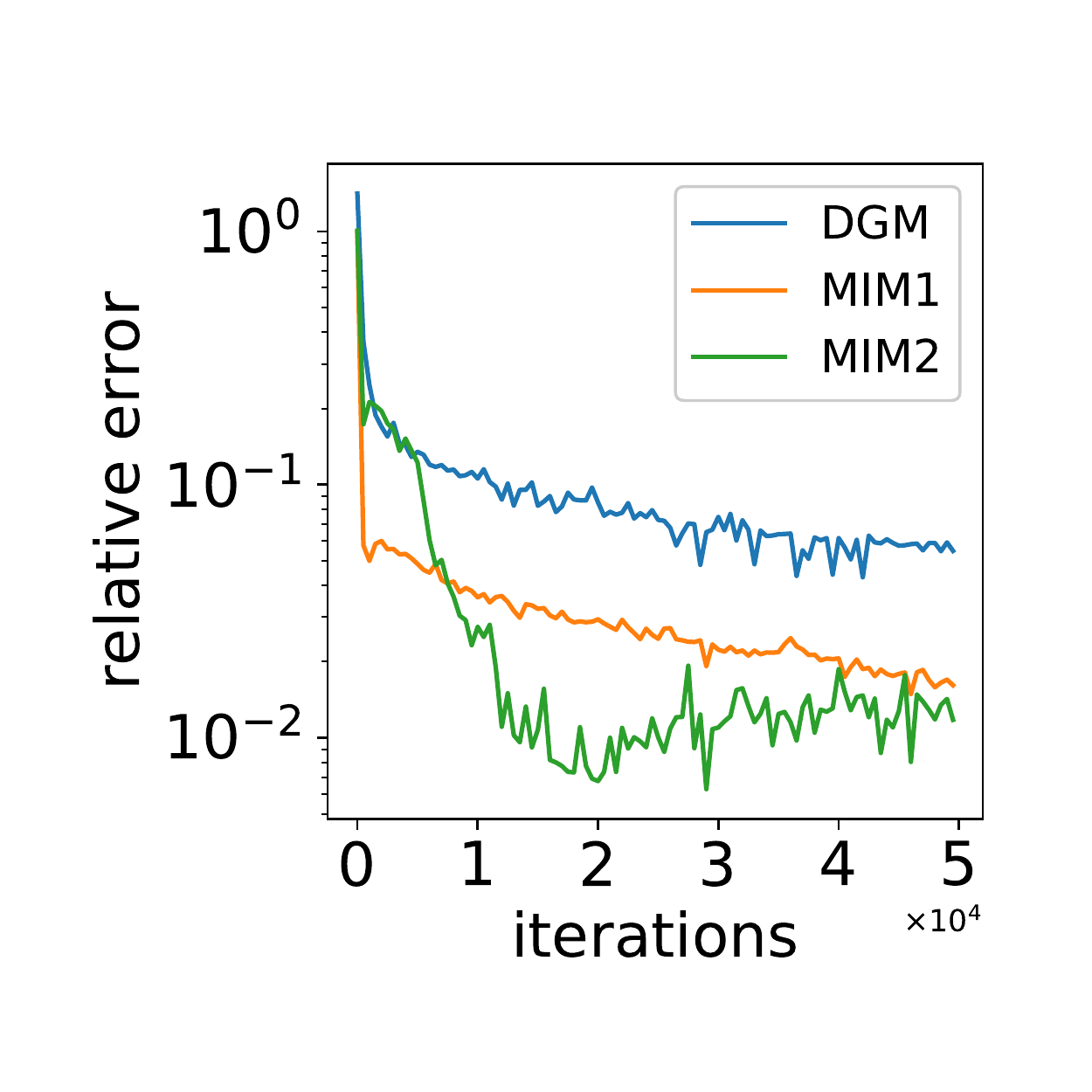}
	}
	\subfigure[$n = 40$, ReCu]{
		\includegraphics[width=0.3\textwidth]{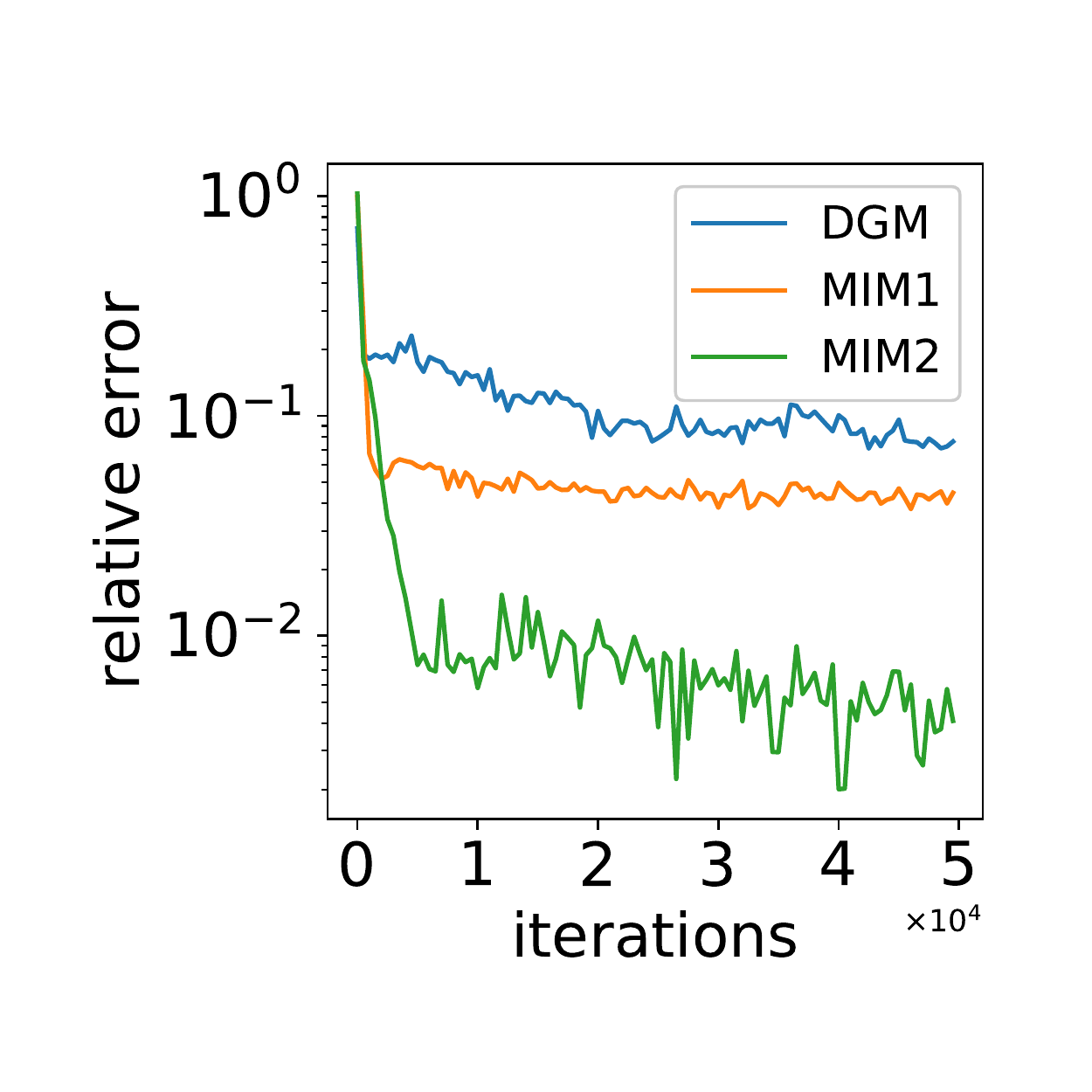}
	}
	\caption{Training processes for wave equation \eqref{eqn:wave} by DGM, MIM1 and MIM2 with the network depth $m = 3$ when $d=3$.}
	\label{fig:wave 3D3L}
\end{figure}

Consider the wave equation
\begin{equation}\label{eqn:wave}
\left\{
\begin{aligned}
& u_{tt} - \Delta u = f & (t,x)\in(0,1)\times\Omega\\
& u(x,t) = 0 & (t,x)\in(0,1)\times \partial \Omega \\
& u(x,0) = 0 & x\in \Omega \\
& u_t(x,0) = 0 & x\in \Omega \\
\end{aligned}\right.
\end{equation}
with the exact solution $u(x,t) = t^2 \Pi_{i=1}^d \sin(\pi x_i)$. It is easy to build the trial DNN that satisfies the initial condition $u(x,0) = 0$ and the boundary condition $u(x,t) = 0$ in both DGM and MIM since both methods have the direct access to the solution $u(x,t)$. However, only MIM is capable of satisfying the other initial condition $u_t(x,0) = 0, \forall x\in\Omega$ since MIM also has the direct access to the derivative by construction
\begin{equation}\label{eqn:wave ic}
\begin{aligned}
&\hat{u}_{\theta} = t\Pi_{i=1}^d(x_i-x_i^2)N_\theta(x), \\
&\hat{v}_{\theta} = t\tilde{N}_\theta(x).
\end{aligned}
\end{equation}
Both $N_\theta(x)$ and $\tilde{N}_\theta(x)$ are one dimensional DNNs.
To compare DGM , MIM1 and MIM2, we add the penalty term to the loss function of DGM with the following form
\begin{equation}
L(u) = \|u_{tt}-\Delta u -f \|^2_{L^2((0,1)\times\Omega)} + \lambda\|u_t \|^2_{L^2((0,1)\times\Omega)}.
\end{equation}
The loss function for MIM1 is
\begin{equation}
\begin{aligned}
L(u,p) = & \|u_{tt}-\nabla\cdot p -f \|^2_{L^2((0,1)\times\Omega)} + \|p-\nabla u \|^2_{L^2((0,1)\times\Omega)} \\
& + \lambda\|u_t \|^2_{L^2((0,1)\times\Omega)}.
\end{aligned}
\end{equation}
The loss function for MIM2 is 
\begin{equation}
\begin{aligned}
L(u,v,p) = & \|v_{t}-\nabla \cdot p -f \|^2_{L^2((0,1)\times\Omega)} + \|p-\nabla u\|^2_{L^2((0,1)\times\Omega)} \\
&+ \|v - u_t \|^2_{L^2((0,1)\times\Omega)}.
\end{aligned}
\end{equation}
We use $\lambda=1$ in DGM and MIM1 for numerical examples. Results of DGM, MIM1, and MIM2 are recorded in Table \ref{tab:wave} with training processes in Figure \ref{fig:wave 3D2L} and Figure \ref{fig:wave 3D3L}. It is clear that MIM2 outperforms DGM and MIM1 in all cases, which is attributed to the enforcement of exact initial conditions in MIM2.
\begin{table}[ht]
	\caption{Relative $L^2$ errors for wave equation \eqref{eqn:wave} with initial conditions. $50000$ sampling points are used in $\Omega$ and the training process ends after $50000$ epochs.}\label{tab:wave}
	\centering
	\begin{tabular}{|c|c|c|c|c|c|}
		\hline
		\multirow{2}*{$d$} & \multirow{2}*{$n$} & \multirow{2}*{$\sigma$}  &\multirow{2}*{Method} & \multicolumn{2}{c|}{$\epsilon$} \\
		\cline{5-6}
		~& ~ & ~ & ~ & $m=2$ & $m=3$\\
		\hline
		\multirow{18}*{2}& \multirow{6}*{10} & \multirow{2}*{ReQu}  
		&DGM &  1.25 e-01 & 7.28 e-02\\
		~& ~ & ~  &MIM1 &  5.20 e-02 & 6.33 e-03\\
		~& ~ & ~  &MIM2 &  7.02 e-02 & 2.90 e-03\\
		\cline{3-6}
		~ & ~ & \multirow{3}*{ReCu}
		&DGM &  1.79 e-02 & 2.39 e-02\\
		~& ~ & ~  &MIM1 &  1.23 e-02 & 3.84 e-03\\
		~& ~ & ~  &MIM2 &  6.89 e-03 & 7.20 e-03\\
		\cline{2-6}
		~ & \multirow{6}*{20} & \multirow{3}*{ReQu}
		&DGM &  4.58 e-02  & 1.68 e-02\\
		~& ~ & ~  &MIM1 &  4.58 e-03  & 2.21 e-03\\
		~& ~ & ~  &MIM2 &  3.71 e-03  & 2.47 e-03\\
		\cline{3-6}
		~ & ~ & \multirow{3}*{ReCu}
		&DGM &  1.87 e-02  & 1.14 e-02\\
		~& ~ & ~  &MIM1 &  1.19 e-03  & 1.13 e-03\\
		~& ~ & ~  &MIM2 &  3.19 e-03  & 2.62 e-03\\
		\cline{2-6}
		~ & \multirow{6}*{40} & \multirow{3}*{ReQu}
		&DGM &  2.77 e-02  & 1.24 e-02\\
		~& ~ & ~  &MIM1 &  1.67 e-03  & 1.42 e-03\\
		~& ~ & ~  &MIM2 &  2.77 e-03  & 2.23 e-03\\
		\cline{3-6}
		~ & ~ & \multirow{3}*{ReCu}
		&DGM &  4.91 e-03  & 3.11 e-03\\
		~& ~ & ~  &MIM &  1.33 e-03  & 1.22 e-03\\
		~& ~ & ~  &MIM &  1.67 e-03  & 1.83 e-03\\
		\hline
		\multirow{18}*{3}& \multirow{6}*{10} & \multirow{3}*{ReQu}
		&DGM &  2.05 e-01 &  1.86 e-01 \\
		~& ~ & ~  &MIM &  2.88 e-02 &  5.85 e-02\\
		~& ~ & ~  &MIM &  1.64 e-02 &  6.21 e-03\\
		\cline{3-6}
		~ & ~ & \multirow{3}*{ReCu}
		&DGM &  1.34 e-01   & 1.30 e-01\\
		~& ~ & ~  &MIM &  5.13 e-02   & 3.17 e-02\\
		~& ~ & ~  &MIM &  2.34 e-02   & 1.47 e-02\\
		\cline{2-6}
		~ & \multirow{6}*{20} & \multirow{3}*{ReQu}
		&DGM     &  1.54 e-01  &  1.01 e-01\\
		~& ~ & ~  &MIM &  4.30 e-02  & 4.03 e-02\\
		~& ~ & ~  &MIM &  1.57 e-02  & 9.32 e-03\\
		\cline{3-6}
		~ & ~ & \multirow{3}*{ReCu}
		&DGM     &  5.66 e-02  &  5.63 e-02\\
		~& ~ & ~  &MIM &  2.23 e-02  & 1.62 e-02\\
		~& ~ & ~  &MIM &  2.02 e-02  & 1.18 e-02\\
		\cline{2-6}
		~ & \multirow{6}*{40} & \multirow{3}*{ReQu}
		&DGM     &  5.98 e-02  &  7.34 e-02  \\
		~& ~ & ~  &MIM &  3.47 e-02  & 4.34 e-03\\
		~& ~ & ~  &MIM &  1.02 e-02  & 4.41 e-03\\
		\cline{3-6}
		~ & ~ & \multirow{3}*{ReCu}
		&DGM     &  1.74 e-02  &  2.15 e-02  \\
		~& ~ & ~  &MIM &  4.01 e-03  & 2.80 e-03\\
		~& ~ & ~  &MIM &  3.27 e-03  & 6.11 e-03\\
		\hline
	\end{tabular}
	
\end{table}
%\begin{figure}
%	\centering
%	\subfigure[$n = 10$, ReQu]{
%		\includegraphics[width=0.3\textwidth]{Wave_2D_w10_l2_requ-eps-converted-to.pdf}
%	}
%	\subfigure[$n = 20$, ReQu]{
%		\includegraphics[width=0.3\textwidth]{Wave_2D_w20_l2_requ-eps-converted-to.pdf}
%	}
%	\subfigure[$n = 40$, ReQu]{
%		\includegraphics[width=0.3\textwidth]{Wave_2D_w40_l2_requ-eps-converted-to.pdf}
%	}
%	\subfigure[$n = 10$, ReCu]{
%	\includegraphics[width=0.3\textwidth]{Wave_2D_w10_l2_recu-eps-converted-to.pdf}
%	}
%	\subfigure[$n = 20$, ReCu]{
%		\includegraphics[width=0.3\textwidth]{Wave_2D_w20_l2_recu-eps-converted-to.pdf}
%	}
%	\subfigure[$n = 40$, ReCu]{
%		\includegraphics[width=0.3\textwidth]{Wave_2D_w40_l2_recu-eps-converted-to.pdf}
%	}
%	\caption{Training processes for wave equation \eqref{eqn:wave} by DGM, MIM1 and MIM2 with the network depth $m = 2$ when $d=2$.}
%	\label{fig:wave 2D2L}
%	\end{figure}
%\begin{figure}
%	\centering
%	\subfigure[$n = 10$, ReQu]{
%		\includegraphics[width=0.3\textwidth]{Wave_2D_w10_l3_requ-eps-converted-to.pdf}
%	}
%	\subfigure[$n = 20$, ReQu]{
%		\includegraphics[width=0.3\textwidth]{Wave_2D_w20_l3_requ-eps-converted-to.pdf}
%	}
%	\subfigure[$n = 40$, ReQu]{
%		\includegraphics[width=0.3\textwidth]{Wave_2D_w40_l3_requ-eps-converted-to.pdf}
%	}
%	\subfigure[$n = 10$, ReCu]{
%		\includegraphics[width=0.3\textwidth]{Wave_2D_w10_l3_recu-eps-converted-to.pdf}
%	}
%	\subfigure[$n = 20$, ReCu]{
%		\includegraphics[width=0.3\textwidth]{Wave_2D_w20_l3_recu-eps-converted-to.pdf}
%	}
%	\subfigure[$n = 40$, ReCu]{
%		\includegraphics[width=0.3\textwidth]{Wave_2D_w40_l3_recu-eps-converted-to.pdf}
%	}
%	\caption{Training processes for wave equation \eqref{eqn:wave} by DGM, MIM1 and MIM2 with the network depth $m = 2$ when $d=3$.}
%	\label{fig:wave 2D3L}
%\end{figure}
\section{Conclusions}
\label{sec:conclusion}

In this work, we propose a systematical strategy to design DNNs that satisfy boundary and initial conditions automatically in the framework of MIM. Since MIM treats both the PDE solution and its derivatives as independent variables, we are able to make DNNs satisfy exact conditions in all cases. Numerous examples are tested to demonstrate the advantages of MIM. Without any penalty term for boundary and initial conditions, MIM does not introduce any modeling error. Therefore, MIM provides better approximations in general, while a deep-learning method with penalty terms typically requires a tuning of penalty parameters in order to produce better results. Note that the penalty term requires an approximation of the $d-1$ dimensional boundary integral which can be prohibitively difficult over a high-dimensional complex domain. Besides, the absence of penalty term facilitates the training process.

\section*{Acknowledgments}
This work is supported in part by the grants National Key R\&D Program of China No. 2018YF645B0204404 and NSFC 21602149 (J.~Chen), and NSFC 11501399 (R.~Du).

\section*{References}
\bibliographystyle{amsplain}
\bibliography{ref}

\end{document}